\documentstyle[12pt]{article}
\begin{document}
\title{Functions of several quaternion variables and
quaternion manifolds.}
\author{S.V. L\"udkovsky}
\date{29 January 2003}
\maketitle
\begin{abstract}
Functions of several quaternion variables are investigated
and integral representation theorems for them are proved.
With the help of them solutions of the ${\tilde \partial }$-equations
are studied. Moreover, quaternion Stein manifolds are defined
and investigated.
\end{abstract}
\section{Introduction}
Superanalysis and its enormous applications in mathematical and theoretical
physics are developing fastly in recent times, but mainly for
supercommutative superalgebras \cite{berez,connes,dewitt,khren}.
For nonsupercommutative superalgebras superanalysis is rather new
\cite{connes,oystaey} and less known, even (super)analysis over
Clifford algebras and, in particular, quaternions is very little
known.
\par Quaternion manifolds appear to be very useful, since
the spin structure for them naturally arise due to the embedding
of the unitary group $U(2)$ of all complex $2\times 2$
unitary matrices into the quaternion skew field $\bf H$:
$U(2)\hookrightarrow \bf H$. Moreover, due to Theorem $4.9$
for each complex manifold $N$ there exists a quaternion manifold
$M$ and a complex holomorphic embedding $\theta : N\hookrightarrow M$.
In view of the isomorphism of the spin group $Spin (4)$
with the direct product of the special unitary groups
$SU(2)\otimes SU(2)$ and the embedding $SU(2)\otimes SU(2)
\hookrightarrow \bf H^2$ each spin manifold $N$ has the embedding
into the corresponding quaternion manifold $M$ of the quaternion
dimension $dim_{\bf H}M=2n$, where $dim_{\bf C}N=2n\in \bf N$
(about complex spin manifolds see \cite{lawmich,moore}).
\par This is possible for manifolds with connecting mappings of
charts being in the class of quaternion holomorphic functions
considered in \cite{luoyst}, but in general this is impossible
in the class of quaternion (right-)superlinearly superdifferentiable
functions, since the latter leads to the additional restrictions
in a complex manifold $N$ associated with the graded $\bf H$ structure
in the tangent space $T_xN$ for each $x\in N$.
The latter manifolds are called in the literature quaternion manifolds
(see, for example, \cite{museya} and references therein),
but practically they are complex manifolds with the additional quaternion
structure in $TN$, since they are modelled on $\bf C^n$ and
their connecting functions of charts are complex holomorphic,
but $\bf H$ is not the algebra over $\bf C$, though there is
the embedding of $\bf C$ into $\bf H$ as the $\bf R$-linear space.
On the other hand, a quaternion manifold $M$ modelled on $\bf H^n$
with quaternion holomorphic connecting functions is not a complex manifold.
Even if $M$ is foliated by the corresponding complex local coordinates
it is like the product of four manifolds $\mbox{ }_{l,m}M$,
$l, m\in \{ 1,2 \} $, with the complex holomorphic structure in
$\mbox{ }_{1,1}M$ and $\mbox{ }_{1,2}M$
and the complex antiholomorphic structure in $\mbox{ }_{2,1}M$ and
$\mbox{ }_{2,2}M$ factorized by the equivalence relation $\cal Y$
such that $z\in \mbox{ }_{1,1}M$ is $\cal Y$-equivalent to $z'\in
\mbox{ }_{2,2}M$ and $\xi \in \mbox{ }_{1,2}M$ is $\cal Y$-equivalent
to $\xi '\in \mbox{ }_{2,1}M$, where $\mbox{ }_{1,1}\phi _a(z)=
{\overline {\mbox{ }_{2,2}\phi _a(z')}}$ and
$\mbox{ }_{1,2}\phi _a(\xi )=-{\overline {\mbox{ }_{2,1}\phi _a(\xi ')}}$,
where $At(M)=\{ (U_a,\phi_a): a\in \Lambda \} $,
$At(\mbox{ }_{l,m}M)=\{ (\mbox{ }_{l,m}U_a,\mbox{ }_{l,m}\phi _a):
a\in \Lambda \} $, $\phi _a: U_a\to \phi _a(U_a)\subset \bf H^n$
is a homeomorphism for each $a$, $\phi _b\circ \phi _a^{-1}$ is
quaternion holomorphic on $\phi _a(U_a\cap U_b)$,
$U_a$ is open in $M$, $\phi _a(U_a)$ is open in $\bf H^n$,
$\bigcup_aU_a=M$.
It is known from the theory of complex manifolds, that
the product $S^{2n+1}\times S^{2m+1}$ of two odd-dimensional
unit real spheres can be supplied with the complex manifold
structure \cite{lnma}. In view of the discussion above
and Theorem $4.9$ the product $S^{2n+1}\times S^{2m+1}
\times S^{2p+1}\times S^{2q+1}$ with $n+m=p+q$ can be supplied
with the quaternion manifold structure, where $n, m, p, q$
are nonnegative integers. In some sense quaternion manifolds
may be related with hyperbolic manifolds.
\par Therefore, in the class of quaternion holomorphic functions
the quantum mechanics on complex manifolds is embeddable
into the quantum mechanics on quaternion manifolds.
Moreover, this is natural, since the Dirac operator can be expressed
through the differentiation by quaternion variables \cite{lawmich}.
Some attempts to spread quantum mechanics from the complex case
into the quaternion case were made in \cite{emch}, but
he had not any mathematical tool of functions of quaternion
variables and quaternion manifolds.
\par In \cite{luoyst} the quaternion line integral $\int_{\gamma }
f(z,{\tilde z})dz$ on the space $C^0_b(U,{\bf H})$
of all bounded continuous functions $f: U\to \bf H$ was investigated,
where $U$ is an open subset in $\bf H$, $\gamma $ is a rectifiable path
in $U$, such that it has some properties like the complex Cauchy
integral: $\int_{\gamma }f(z,{\tilde z})dz=\int_{\gamma _1}
f(z,{\tilde z})dz+\int_{\gamma _2}f(z,{\tilde z})dz$ for each
$\gamma =\gamma _1\cup \gamma _2$ with $\gamma _1(1)=\gamma _2(0)$,
$\gamma _1, \gamma _2: [0,1]\to U$, $\int_{\gamma }(af(z,{\tilde z})
+bg(z,{\tilde z})dz=a\int_{\gamma }f(z,{\tilde z})dz+b\int_{\gamma }
g(z,{\tilde z})dz$ for each $a, b\in \bf H$, $f, g \in C^0_b(U,{\bf H})$,
$(\partial (\int_{\gamma ,\eta =\gamma (1)}f(z,{\tilde z})dz)/
\partial \eta ).I=f(\eta ,{\tilde \eta })$, though the quaternion line
integral is not $\bf H$-linear. Then there quaternion holomorphic
functions were introduced and investigated with the specific
definition of the superdifferentiability (in general $\bf H$-nonlinear).
There were found many quaternionic
features in comparison with the complex case
and the theory of quaternion holomorphic functions can not be reduced
to the theory of complex holomorphic functions.
\par In this paper functions of several quaternion variables
are investigated (see \S \S 2 and 3) and their theory is applied
to the definition and the studies of quaternion analogs
of complex Stein manifolds (see \S 4).
In this article theorems about integral representations
of quaternion functions are proved.
Among them there are the quaternion analogs of the Cauchy-Green,
Martinelli-Bochner and Leray formulas, but they are quite different
from the complex case, since quaternion integrals
are noncommutative and quaternion differential forms does not
satisfy the same properties as complex differential forms.
They are applied to solve the $\tilde \partial $-equations.
This is important not only for the theory of quaternion
functions, but also for investigations of quaternion manifolds.
Then this can be used for developments of sheaves, quantum sheaves
and quantum field theory on quaternion manifolds. Naturally
free loop spaces and transformation groups of quaternion manifolds
such as groups of diffeomorphisms, groups of geometric loops,
groups of hoops can be further investigated continuing previous
works of the author on these spaces and groups
\cite{lulgcm,lulsqm,lustptg,lupm}, which is interesting
for theory of groups and their representations and also for their
usages in quantum field theory, quantum gravity, superstrings, etc.
\section{Differentiable functions of several quaternion variables}
\par {\bf 2.1. Theorem.} {\it Let $U$ be an open subset in $\bf H$
with a $C^1$-boundary $\partial U$ $U$-homotopic with a product
$\gamma _1\times \gamma _2\times \gamma _3$, where $\gamma _j(s)=
a_j+r_j\exp (2\pi M_js)$, $M_j\in \bf H_i$, $|M_j|=1$, $s\in [0,1]$,
$\gamma _j([0,1])\subset U$, $0<r_j<\infty $, $j=1,2,3$, where
$M_j$ are linearly independent over $\bf R$.
Let also $f: cl (U)\to \bf H$ be a continuous function on $cl (U)$
such that $(\partial f(z)/\partial {\tilde z})$ is defined in the sense
of distributions in $U$ is continuous in $U$ and has a continuous
extension on $cl (U)$, where $U$ and $\gamma _j$ for each
$j$ satisfy conditions of Theorem 3.9 \cite{luoyst}. Then
$$(2.1) \quad f(z)=$$
$$(2\pi )^{-3}\int_{\gamma _3}\int_{\gamma _2}
\int_{\gamma _1}
f(\zeta _1)(\partial _{\zeta _1}Ln(\zeta _1-\zeta _2))M_1^{-1}
(\partial _{\zeta _2}Ln (\zeta _2-\zeta _3))M_2^{-1}(\partial _{\zeta _3}
Ln (\zeta _3-z))M_3^{-1}$$   $$-(2\pi )^{-3}
\int_U \{ (\partial {\hat f}(\zeta _1)/
\partial {\tilde \zeta }_1).d{\tilde \zeta }_1\wedge
\partial _{\zeta _1}Ln(\zeta _1-\zeta _2))M_1^{-1}
(\partial _{\zeta _2}Ln (\zeta _2-\zeta _3))M_2^{-1}(\partial _{\zeta _3}
Ln (\zeta _3-z))M_3^{-1} \} . $$ }
\par {\bf Proof.} We have the identities
$d_{\zeta } [{\hat f}(\zeta ). \partial _{\zeta } Ln (\zeta -z)]=$
$\{ (\partial {\hat f}(\zeta )/\partial {\tilde \zeta }).d{\tilde \zeta } \}
\wedge \partial _{\zeta }Ln (\zeta -z)$ $+\{
(\partial {\hat f} (\zeta )/\partial {\zeta }).d\zeta \}
\wedge \partial _{\zeta } Ln (\zeta -z) $ and
$d_{\zeta }d_{\zeta }Ln (\zeta -z)|_{\zeta \in \gamma }=0$
for $\zeta $ varying along a path $\gamma $,
where for short $f(z)=f(z,\tilde z)$, since there is the bijection
of $z$ with $\tilde z$ on $\bf H$,
there exist quaternion-valued functions $g(\zeta ,z)$ and
$h(\zeta ,z)$ such that $\partial _{\zeta }Ln (\zeta -z)=
g(\zeta ,z)d\zeta =(d\zeta )h(\zeta ,z)$ (see also
\S \S 2.1 and 2.6 \cite{luoyst}).
As in \cite{luoyst} ${\hat f}(z,{\tilde z}):=\partial g(z,{\tilde z})/
\partial z$, where $g(z,{\tilde z})$ is a quaternion-valued
function such that $(\partial g(z,{\tilde z})/\partial z).I=
f(z,{\tilde z})$. Since $\zeta _1$ varies along the path
$\gamma _1$, then $d\zeta _1\wedge d\zeta _1|_{\zeta _1\in \gamma _1}=0$.
Consider $z\in U$ and $\epsilon >0$ such that the torus
${\bf T}(z,\epsilon ,{\bf H})$ is contained in $U$, where
$\partial {\bf T} (z,\epsilon ,{\bf H}) =\psi _3\times \psi _2
\times \psi _1$, $\psi _j$ are of the same form as $\gamma _j$ but with
$z$ instead of $a_j$ and with $r_j=\epsilon $.
Applying Stokes formula for regions in $\bf R^4$ and componentwise
to $\bf H$-valued differential forms we get
$$\int_{\partial U}\omega -\int_{\partial {\bf T}(z,\epsilon ,{\bf H})}
\omega =
\int_{U\setminus {\bf T}(z,\epsilon ,{\bf H})}dw, \mbox{ where}$$
$$w={\hat f}(\zeta _1).(\partial _{\zeta _1}Ln(\zeta _1-\zeta _2))M_1^{-1}
(\partial _{\zeta _2}Ln (\zeta _2-\zeta _3))M_2^{-1}(\partial _{\zeta _3}
Ln (\zeta _3-z))M_3^{-1}.$$ In view of Theorem 3.9 \cite{luoyst}
we have that $\lim_{\epsilon \to 0, \epsilon >0}
(2\pi )^{-3}\int_{\gamma _3}\int_{\gamma _2} \int_{\gamma _1}w=f(z)$
and $\lim_{\epsilon \to 0, \epsilon >0}
\int_{U\setminus {\bf T}(z,\epsilon ,{\bf H})}d\omega =\int_Ud\omega .$
From this formula $(2.1)$ follows.
\par {\bf 2.1.1. Remark.} Formula $(2.1)$ is the quaternion analog
of the (complex) Cauchy-Green formula.
Since in the sence of distributions $\partial {\hat f}/\partial {\tilde z}=
\partial (\partial g/\partial {\tilde z})/\partial z$,
then from $\partial {\hat f}/\partial {\tilde z}=0$ it follows, that
$\partial {\hat f}.I/\partial {\tilde z}=\partial f/\partial {\tilde z}=0$.
If $\partial f/\partial {\tilde z}=0$, then $g$ can be chosen such that
$\partial g/\partial {\tilde z}=0$ \cite{luoyst}.
Therefore, from Formula $(2.1)$ it follows, that $f$ is quaternion
holomorphic in $U$ if and only if $\partial f/\partial {\tilde z}=0$
in $U$.
\par {\bf 2.2. Remark.} Instead of curves $\gamma $ of Theorem
$2.1$ above or Theorems $3.9, 3.24$ \cite{luoyst} it is possible
to consider their natural generalization
$\gamma (\theta )+z_0=z_0+r(\theta ) \exp (2\pi S(\theta ))$,
where $r(\theta )$ and $S(\theta )$ are continuous functions
of finite total variations, $\theta \in [0,1]\subset \bf R$,
$r(\theta )\ge 0$, $S(\theta )\in \bf H$.
Therefore, $\gamma $ is a rectifiable path.
If $S(0)=S(1)\quad mod({\bf S_i})$ and $r(0)=r(1)$, then $\gamma $
is a closed path (loop): $\quad \gamma (0)=\gamma (1)$, where
${\bf S_i}:=\{ z\in {\bf H_i}: |z|=1 \} $,
${\bf H_i}:=\{ z\in {\bf H}: z+{\tilde z}=0 \} $.
Consider $S$ absolutely continuous such that there exists
$T\in L^1([0,1],{\bf H})$ for which $S(\theta )=
S(0)+\int_0^{\theta }T(\tau )d\tau $
(see Satz 2 and 3 (Lebesgue) in \S 6.4 \cite{kolmfom})
and let $r(\theta )>0$ for each $\theta \in [0,1]$.
Evidently, $Mn:=S(1)-S(0)=\int_0^1T(\tau )d\tau $ is invariant relative
to reparatmetrizations $\phi \in Diff^1_+([0,1])$ of diffeomorphisms
of $[0,1]$ preserving the orientation, $n$ is a real number,
$M\in {\bf S_i}$. Then $\Delta Arg (\gamma )
:=Arg (\gamma )|_0^1= 2\pi \int_0^1T(\tau )d\tau $
(see also Formula $(3.7)$ and \S 3.8 \cite{luoyst}).
In view of Theorem $3.8$ \cite{luoyst} for each loop
$\gamma :$ $\quad \Delta Arg (\gamma )\in \bf ZS_i$.
For each $\epsilon >0$ for the total variation there is
the equality $V(\gamma \epsilon )=V(\gamma )\epsilon $.
Since $\gamma ([0,1])$ is a compact subset in $\bf H$,
then there exists $r_m:=\sup _{\theta \in [0,1]}r(\theta )<\infty $.
Hence $z_0+(\gamma \epsilon )([0,1])\subset
B({\bf H},z_0,r_m\epsilon )$.
\par Therefore, Theorems $3.9, 3.24, 3.30$ and Formulas
$(3.9, 3.9')$ \cite{luoyst} and Theorem $2.1$ above are true for such
paths $\gamma $ also and Formula $(3.9)$ \cite{luoyst} takes the form
$$(2.2)\quad f(z)=(2\pi n)^{-1}(\int_{\psi }f(\zeta )
(\zeta -z)^{-1}d\zeta ){\tilde M},$$
where $0\ne n\in \bf Z$ for a closed path $\gamma $, $M\in \bf S_i$,
Formula $(2.2)$ generalizes Formula $(3.9)$, when $|n|>1$.
When ${\hat I}n(0,\gamma )=0$, then
$(\int_{\psi }f(\zeta )(\zeta -z)^{-1}d\zeta )=0$
(see also \S 3.23 \cite{luoyst}).
\par {\bf 2.3. Theorem.} {\it Let $U$ be a bounded open subset
in $\bf H$ and $f: U\to \bf H$ be a bounded continuous function.
Then there exists a continuous function $u(z)$ which
is a solution of the equation
$$(1)\quad (\partial u(z)/\partial {\tilde z})={\hat f}$$
in $U$, in particular, $(\partial u(z)/\partial {\tilde z}).I=f(z)$.}
\par {\bf Proof.} Take quaternion one-forms $d\zeta _l$ expressible through
$d\zeta $ as $\sum_{j=1}^{k(l)}P_{j,1,l}d\zeta P_{j,2,l}$ with
fixed nonzero quaternions $P_{j,q,l}$, where $l=1,2,3,4$, $k(l)\in \bf N$,
$\zeta \in U$.
Choose them satisfying conditions $d\zeta _4\wedge \nu =\xi (v,w,x,y)
dv\wedge dw\wedge dx\wedge dy$, $d{\tilde \zeta _4}\wedge \nu =0$,
where $\nu =(\partial _{\zeta _1} Ln (\zeta _1-\zeta _2))M_1^{-1}
(\partial _{\zeta _2} Ln (\zeta _2-\zeta _3))M_2^{-1}
(\partial _{\zeta _3} Ln (\zeta _3-z))M_3^{-1}$
as in \S 2.1, $z=vI+wJ+xK+yL$, $I$, $J$, $K$ and $L$
are Pauli-matrices, $v$, $w$, $x$ and $y\in \bf R$,
$\xi : U\to \bf \hat H$ is a function nonzero and finite
almost everywhere on $U$ relative to the Lebesgue measure.
Then there exist $d\zeta _4$ and $\nu $ such that the continuous function
$$(2)\quad u(z):=-(2\pi )^{-3}\int_U ({\hat f}
(\zeta _1).d{\tilde \zeta }_4 ) \wedge \nu $$
is a solution of equation $(1)$. To demonstrate this
take closed curves $\gamma _j$
in $U$ as in \S 2.1 and \S 2.2, for example, such that
$\zeta _j\in \gamma _j$ satisfy conditions:
$(\zeta _1-\zeta _2)=\eta _1$, $(\zeta _2-\zeta _3)=\eta _3$,
$\zeta _3=r\eta _3^T$ with $0<r<1$, where $\eta _1={{t\quad 0}\choose
{0\quad {\bar t}}}$, $\eta _2={\tilde \eta }_3$,
$\eta _3={{0\quad u}\choose {-{\bar u}\quad 0}}$,
${\tilde \eta }_1=\eta _4$,
$\eta _4={{{\bar t}\quad 0}\choose {0\quad t}}$, where
$t$ and $u\in \bf C$, $\zeta ={{t\quad u}\choose
{{-\bar u}\quad {\bar t}}}\in U$,
$a^T$ denotes the transposed matrix of a matrix $a$.
Hence $d\eta _1\wedge d\eta _1=0$,
$d\eta _2\wedge d\eta _2^T=0$, $d\eta _3\wedge d\eta _3^T=0$,
$d\eta _4\wedge d\eta _4=0$,
\par $(i)$ $\eta _1^kd\eta _1=(d\eta _1)\eta _1^k$ and
$\eta _3^kd\eta _3=[(d\eta _3)\eta _3^k]^{\tilde .}$ \\
for $k=1$ and $k=-1$. These variables are expressible as
$\zeta _l=\sum_{j=1}^{k(l)}P_{j,1,l}\zeta P_{j,2,l}$
(see \S \S 3.7 and 3.28 \cite{luoyst}).
Therefore, there exists a quaternion variable $\zeta _4$ expressible through
$\zeta $ as above such that
\par $(ii)$ $d{\tilde \zeta }_4\wedge \nu =\xi (v,w,x,y)
dv\wedge dw\wedge dx\wedge dy$, $d\zeta _4 \wedge \nu =0$.
Therefore, there exists a subgroup of the group of all quaternion
holomorphic diffeomorphisms of $U$ preserving Conditions $(ii)$
and the construction given above has natural generalizations.
\par On the other hand, $\zeta $ is expressible through
$\zeta _1$,...,$\zeta _4$ as $\zeta =\sum_lP_{1,l}\zeta P_{2,l}$,
where $P_{1,l},...,P_{2,l}$ are quaternion constants.
Let at first $f$ be continuously differentiable in $U$.
Each $\zeta _j$ is expressible in the form $\zeta _j=
\sum_l\mbox{ }^lb_jS_l$, where $\mbox{ }^lb_j\in \bf R$
are real variables, $S_l\in \{ I,J,K,L \} $, hence differentials
$(\partial f/\partial \zeta _j).d\zeta _j=\sum_l \{ (\partial
f/\partial z).S_ld\mbox{ }^lb_j$ $+(\partial f/\partial {\tilde z}).
{\tilde S_l}d\mbox{ }^lb_j \} $ are defined.
Consider a fixed $z_0\in U$. We take a $C^{\infty }$-function
$\chi $ on $\bf H$ such that $\chi =1$ in a neighbourhood $V$ of $z_0$,
$V\subset U$, $\chi =0$ in a neighbourhood of ${\bf H}\setminus U$.
Then $u=u_1+u_2,$ where 
$$u_1(z):=-(2\pi )^{-3}\int_U \chi (\zeta _1)
({\hat f}(\zeta _1).d{\tilde \zeta }_4)\wedge \nu ,$$
$$u_2(z):=-(2\pi )^{-3}\int_U (1-\chi (\zeta _1))
({\hat f}(\zeta _1).d{\tilde \zeta }_4)\wedge \nu .$$
Then
$$u(z):=-(2\pi )^{-3}\int_{\bf H}\chi (\zeta _1+z)
({\hat f}(\zeta _1+z).d{\tilde \zeta }_4)\wedge \psi ,\mbox{ where}$$
$$\psi :=(\partial _{\zeta _1} Ln (\zeta _1-\zeta _2))M_1^{-1}
(\partial _{\zeta _2} Ln (\zeta _2-\zeta _3))M_2^{-1}
(\partial _{\zeta _3} Ln (\zeta _3)M_3^{-1}.$$
Since $\partial _{\zeta _4}
\{ [\chi (\zeta _1+z){\hat f}(\zeta _1+z).]\wedge \psi \} =0$
and $(\partial /\partial {\tilde \zeta _4})
\{ [\chi (\zeta _1+z){\hat f}(\zeta _1+z).]\wedge \psi \}.d{\tilde \zeta }_4$
$=\partial _{{\tilde \zeta }_4}
\{ [\chi (\zeta _1+z){\hat f}(\zeta _1+z).]\wedge \psi \},$
then due to Equations $(i,ii)$ \\
$(\partial u(z)/\partial {\tilde z})=
-(2\pi )^{-3}\int_{\bf H}\partial _{{\tilde \zeta }_4}
\{ [\chi (\zeta _1+z){\hat f}(\zeta _1+z)]\wedge \psi \} $.  \\
In view of Theorem 2.1 applied to ${\hat f}.S$ for each
$S\in \{ I,J,K,L \} $ we have
$(\partial u_1/\partial {\tilde z})=\hat f$ in $V$,
consequently, $(\partial u/\partial {\tilde z})=\hat f$ in a neighbourhood
of $z_0$.
Taking a sequence $f^n$ of continuously differentiable functions uniformly
converging to $f$ on $U$ we get the corresponding $u^n$ such that
in the sence of distributions $(\partial u/\partial {\tilde z})=
\lim_{n\to \infty } (\partial u^n/\partial {\tilde z})
=\lim_n{\hat f}^n=\hat f$.
\par {\bf 2.4. Theorem.} {\it Let $U$ be an open subset in $\bf H^n$.
Then for every compact subset $K$ in $U$ and every multi-order
$k=(k_1,...,k_n)$, there exists a constant $C>0$ such that
$$\max_{z\in K} |\partial ^kf(z)|\le C\int_U|f(z)|d\sigma _{4n}$$
for each quaternion holomorphic function $f$, where $d\sigma _{4n}$
is the Lebesgue measure in $\bf H^n$.}
\par {\bf 2.5. Corollary.} {\it Let $U$ be an open subset in $\bf H^n$
and let $f_l$ be a sequence of quaternion holomorphic functions
in $U$ which is uniformly bounded on every compact subset
of $U$. Then there is a subsequence $f_{k_j}$ converging
uniformly on every compact subset of $U$ to a limit in
$C^{\omega }_z(U,{\bf H})$.}
\par Proofs of Theorem 2.4 and Corollary 2.5 follow from
Theorem 2.1 above and Theorem 3.9 \cite{luoyst}
analogously to Theorem 1.1.13 and Corollary 1.1.14 \cite{henlei}.
\par {\bf 2.6. Definitions.} Let $U$ be an open subset in $\bf H^n$
and $f: U\to \bf H^m$ be a quaternion holomorphic function,
then the matrix: $J_f(z):=(\partial f_j(z)/\partial z_k)$
is called the quaternion Jacobi matrix, where $j=1,...,m$,
$k=1,...,n$. To this quaternion operator matrix there corresponds
a real $(4m)\times (4n)$-matrix. Denote by $rank_{\bf R}(J_f(z))$
a rank of a real matrix corresponding to $J_f(z)$.
Then $f$ is called regular at $z\in U$,
if $rank_{\bf R}(J_f(z))=4\min (n,m)$.
If $U$ and $V$ are two open subsets in $\bf H^n$, then
a bijective surjective mapping $f: U\to V$ is called
quaternion biholomorphic if
$f$ and $f^{-1}: V\to U$ are quaternion holomorphic.
\par {\bf 2.7. Proposition.} {\it Let $U$ and $V$ be open subsets
in $\bf H^n$ and $\bf H^m$ respectively. If $f: U\to \bf H^m$
and $g: V\to \bf H^k$ are quaternion holomorphic functions
such that $f(U)\subset V$, then $g\circ f: U\to \bf H^k$ is
quaternion holomorphic and $J_{g\circ f}(z)=J_g(f(z)).(J_f(z).h)$
for each $h\in \bf H^n$.}
\par {\bf Proof.} In view of Definition 2.1 and Theorem 3.10
\cite{luoyst}
$(\partial g_j(f(z))/\partial z_l).\zeta =
\sum_{s=1}^m\sum_{l=1}^k (\partial g_j(\xi )/ \partial \xi _s)|_{\xi =f(z)}.
(\partial f_s(z)/\partial z_l).h_l$, where $h=(h_1,...,h_n)$, $h_l\in \bf H$
for each $l=1,...,n$, since $f(U)\subset V$ and this is evident
for quaternion polynomial functions and hence for locally converging
series of quaternion holomorphic functions.
\par {\bf 2.8. Proposition.} {\it Let $U$ be a neighbourhood of
$z\in \bf H^n$ and let $f: U\to \bf H^n$ be a quaternion holomorphic
function. Then $f$ is quaternion biholomorphic in some neighbourhood
of $z$ if and only if $rank_{\bf R}J_f(z)=4n$.}
\par {\bf Proof.} From Proposition 2.7 it follows, that
the condition $rank_{\bf R}J_f(z)=4n$ is necessary.
In view of Definition 2.1, Theorem 3.10 and Note 3.11 \cite{luoyst}
an incerement of $f$ can be written
in the form $f(z+\zeta )=f(z)+J_f(z).\zeta + O(|\zeta |^2)$
for each $\zeta \in \bf H^n$ such that $z+\zeta \in U$.
As in the proof of Theorem 1.1.18 \cite{henlei} we get, that
there exists a quaternion holomorphic function $h$ on an open
neighbourhood $W$ of $z$ in $U$ satisfying the condition
$|g(z+\zeta )|\le C|\zeta |^2$, $W\supset B(z,2\epsilon ,{\bf H^n})$,
$0<\epsilon <(2C)^{-1}$, where $C$ is a positive constant
and $h$ is given by the series
$h=\sum_{k=1}^{\infty }g_k$, where $g_{k+1}=g\circ g_k$
for each $k\in \bf N$ and $g_1:=g$, $g:=id-f$, since
for each $\eta \in W$ there exists $r>0$ such that
$B(\eta ,r,{\bf H})\subset U$ and the series
for $h$ is convergent on $B(\eta ,r,{\bf H})$
with $h(B(z,\epsilon ,{\bf H^n}))\subset B(z,2\epsilon ,{\bf H^n})$.
The operator $J_f(\eta )$ is continuous by $\eta $ on $U$, hence
there exists a neighbourhood $V$ of $z$ such that $rank_{\bf R}J_f(z)$
is equal to $4n$ on it, hence $f(V)$ is open in $\bf H^n$. 
From $(id+h)\circ f=f\circ (id+h)=id$ on $B(z,\epsilon ,{\bf H^n})$
it follows, that $f$ is quaternion holomorphic on a neighbourhood of $z$.
\par {\bf 2.9. Corollary.} {\it Let $X$ be a subset in $\bf H^n$
and $k\in \{ 1,2,...,n-1 \} $, then the following conditions are equivalent:
\par $(i)$ for each $\zeta \in X$ there exists a quaternion
biholomorphic map $f=(f_1,...,f_n)$ in some neighbourhood
$U$ of $\zeta $ such that $rank_{\bf R}f=4n$ on $U$ and
$X\cap U= \{ z\in U: f_{k+1}(z)=0,...,f_n(z)=0 \} $;
\par $(ii)$ for each $\zeta \in X$ there exists a neighbourhood
$V$ of $\zeta $ and a regular quaternion holomorphic map
$g: V\to \bf H^{n-k}$ such that $X\cap V= \{ z\in V: g(z)=0 \} $.}
\par {\bf Proof} is analogous to the proof of Corollary 1.1.19
\cite{henlei} and with a consideration of the determinant
function of the real $4n\times 4n$ matrix $J_{f,{\bf R}}(z)$
corresponding to the quaternion operator $J_f(z)$, since
$det J_{f,{\bf R}}(z)\ne 0$ if and only if $rank_{\bf R}f(z)=4n$.
\par {\bf 2.10. Definitions.} Let $U$ be an open subset in $\bf H^n$.
A subset $X$ in $U$ is called a quaternion submanifold of
$\bf H^n$ if the equivalent conditions of Corollary 2.9 are
satisfied. If in addition $X$ is a closed subset in $U$,
then $X$ is called a closed quaternion submanifold of $U$.
This definition is the particular case of the following general
definition.
\par A quaternion holomorphic manifold of quaternion dimension $n$ is
a real $4n$-dimensional $C^{\infty }$-manifold $X$ together
with a family $ \{ (U_j,\phi _j): j\in \Lambda \} $
of charts such that
\par $(i)$ each $U_j$ is an open subset in $X$ and
$\bigcup_{j\in \Lambda }U_j=X$, where $\Lambda $ is a set;
\par $(ii)$ for each $j\in \Lambda $ a mapping
$\phi _j: U_j\to V_j$ is a homeomorphism on an open subset
$V_j$ in $\bf H^n$;
\par $(iii)$ for each $j, l\in \Lambda $ a connection mapping
$\phi _j\circ \phi _l^{-1}$ is a quaternion biholomorphic map
(see \S 2.6) from $\phi _l(U_j\cap U_l)$ onto $\phi _j(U_j\cap U_l)$.
Such system is called a quaternion holomorphic atlas
$At(X):= \{ (U_j,\phi _j): j\in \Lambda \} $.
Each chart $(U_j,\phi _j)$ provides a system of quaternion
holomorphic coordinates induced from $\bf H^n$.
For short we shall write quaternion manifold instead of
quaternion holomorphic manifold and quaternion atlas instead
of quaternion holomorphic atlas if other will not be specified.
\par For two quaternion manifolds $X$ and $Y$
with atlases $At(X):= \{ (U_j,\phi _j): j\in \Lambda _X \} $
and $At(Y):= \{ (W_l,\psi _l): l\in \Lambda _Y \} $
a function $f: X\to Y$ is called quaternion holomorphic if
$\psi _l\circ f\circ \phi _j^{-1}$ is quaternion holomorphic
on $\phi _j(U_j\cap f^{-1}(W_l))$.
If $f: X\to Y$ is a quternion biholomorphic epimorphism,
then $X$ and $Y$ are called quaternion biholomorphically
equivalent.
\par A subset $Z$ of a quaternion manifold $X$ is called
a quaternion submanifold, if $\phi _j(U_j\cap Z)$
is a quaternion submanifold in $\bf H^n$ for each chart
$(U_j,\phi _j)$. If additionally $Z$ is closed in $X$, then
$Z$ is called a closed quaternion submanifold.
\par {\bf 2.11. Theorem.} {\it Let $n\ge 2$, $f_1,...,f_n
\in C^1_{0,(z,\tilde z)}({\bf H^n},{\bf H})$ be a family of continuously
quaternion $(z,\tilde z)$-differentiable functions
satisfying compatibility conditions:
$$(i)\quad \partial f_j/\partial {\tilde z}_k=
\partial f_k/\partial {\tilde z}_j\mbox{ for each }
j, k=1,...,n,$$
where in $C^1_{0,(z,\tilde z)}({\bf H^n},{\bf H})$ is the subspace
of $C^1_{(z,\tilde z)}({\bf H^n},{\bf H})$ of functions with compact support.
Then there exists $u\in C^1_{0,(z,\tilde z)}({\bf H^n},{\bf H})$
satisfying the following $\tilde \partial $-equation:
$$(ii)\quad \partial u/\partial {\tilde z}_j={\hat f}_j,
\quad j=1,...,n; $$
in particular, $(\partial u/\partial {\tilde z}_j).I=f_j.$}
\par {\bf Proof.} We put
$$(iii)\quad u(z):=-(2\pi )^{-3}\int_{\bf H}
[({\hat f}_1(\zeta _1,z_2,...,z_n).d{\tilde \zeta }_4)
\wedge \eta ],\mbox{ where}$$
$$\eta :=(\partial _{\zeta _1} Ln (\zeta _1-\zeta _2))M_1^{-1}
(\partial _{\zeta _2} Ln (\zeta _2-\zeta _3))M_2^{-1}
(\partial _{\zeta _3} Ln (\zeta _3-z)M_3^{-1}$$
(see \S 2.3). By changing of variables we get
$$u(z):=-(2\pi )^{-3}\int_{\bf H}[({\hat f}_1(z_1+\zeta _1,z_2,...,z_n).
d{\tilde \zeta }_4)\wedge \psi ],\mbox{ where}$$
$$\psi :=(\partial _{\zeta _1} Ln (\zeta _1-\zeta _2))M_1^{-1}
(\partial _{\zeta _2} Ln (\zeta _2-\zeta _3))M_2^{-1}
(\partial _{\zeta _3} Ln (\zeta _3)M_3^{-1}.$$
Therefore, $u\in C^1_{(z,\tilde z)}({\bf H^n},{\bf H})$.
Due to Theorem 2.3  $\partial u/\partial {\tilde z}_1={\hat f}_1$
in $\bf H^n$. In view of Theorem 2.1 and the condition
$\partial f_1/\partial {\tilde z}_k=\partial f_k/\partial {\tilde z}_1$
the following equality is satisfied
$${\hat f}_k(z)=-(2\pi )^{-3}\int_{\bf H} \{
[\partial f_k(\zeta _1,z_2,...,z_n)/\partial {\tilde \zeta }_1].
d{\tilde \zeta }_4 \} \wedge \psi ,$$
hence $\partial u/\partial {\tilde z}_k={\hat f}_k$ for $k=2,...,n$,
that is, $u$ satisfies equations $(ii)$.
From this it follows, that $u$ is quaternion holomorphic
in ${\bf H^n}\setminus supp (f_1)\cup ... \cup supp (f_n)$.
In view of formula $(iii)$ it follows, that there exists
$0<r<\infty $ such that
\par $(iv) \quad u(z)=0$ for each $z\in \bf H^n$
with $|z_2|+...+|z_n|>r$. From $\partial u/\partial {\tilde z}_1={\hat f}_1$
it follows, that $\partial u/\partial {\tilde z}_1=0$
in ${\bf H^n}\setminus supp (f_1)$. Consequently,
there exists $0<R<\infty $ such that
$u$ may differ from $0$ on ${\bf H^n}\setminus B({\bf H^n},0,R)$
only on a quaternion constant (see
Theorem 3.28 and Note 3.11 in \cite{luoyst}).
Together with $(iv)$ this gives,
that $u(z)=0$ on ${\bf H^n}\setminus B({\bf H^n},0,\max (R,r)).$
\par {\bf 2.12. Theorem.} {\it Let $U$ be an open subset in $\bf H^n$,
where $n\ge 2$. Suppose $K$ is a compact subset in $U$ such that
$U\setminus K$ is connected.
Then for every quaternion holomorphic function $h$ on $U\setminus K$
there exists a function $H$ quaternion holomorphic in $U$ such that
$H=h$ in $U\setminus K$.}
\par {\bf Proof.} Take any infinite $(z,\tilde z)$-differentiable
function $\chi $ on $U$ with compact support such that
$\chi |_V=1$ on some (open) neighbourhood $V$ of $K$.
Then consider a family of functions $f_j$ such that
${\hat f}_j(z,\tilde z).S=
-\{ (\partial \chi /\partial {\tilde z}).S \} h$ in $U\setminus K$
and $f_j=0$ outside $U\setminus K$
for each $S\in \{ I,J,K,L \} $, where $j=1,...,n$,
$f_j(z)={\hat f}_j(z).I$.
Therefore, conditions of Theorem 2.11 are satisfied
and it gives a function $u\in C^1_{0,(z,\tilde z)}({\bf H^n},{\bf H})$
such that $\partial u/\partial {\tilde z}_j={\hat f}_j$ for each
$j=1,...,n$. A desired function $H$ can be defined
by the formula $H:=(1-\chi )h-u$ such that $H$ is quaternion
holomorphic in $U$. Since $\chi $ has a compact support, then there
exists an unbounded connected subset $W$ in ${\bf H^n}\setminus
supp (\chi )$. Therefore, $u|_W=0$, consequently,
$H|_{U\cap W}=h|_{U\cap W}$. From $(U\setminus K)\cap W\ne \emptyset $
and connectedness of $U\setminus K$ it follows, that
$H|_{U\setminus K}=h|_{U\setminus K}$.
\par {\bf 2.13. Remark.} In the particular case of a singleton
$K=\{ z \} $ Theorem 2.12 gives nonexistence of isolated singularities,
that is, each quaternion holomorphic function in $U\setminus \{ z \} $
for $U$ open in $\bf H^n$ with $n\ge 2$ can be quaternion
holomorphically extended to $z$. Theorem 2.12 is the quaternion
analog of the Hartog's theorem for $\bf C^n$.
\par {\bf 2.14. Corollary.} {\it Let $U$ be an open
connected subset in $\bf H^n$
and $n\ge 2$. Suppose that $f$ is a right superlinearly superdifferentiable
function $f: U\to \bf H$ and $N(f):=\{ z\in U: f(z)=0 \} $, then
\par $(i)$ $U\setminus N(f)$ is connected,
\par $(ii)$ $N(f)$ is not compact.}
\par {\bf Proof.} $(i)$. We have $f={{f_{1,1}\quad f_{1,2}}\choose
{-{\bar f}_{1,2}\quad {\bar f}_{1,1}}}$, where
$f_{1,1}(t,u)$ is holomorphic by $t$ and antiholomorphic by $u$,
$f_{1,2}(t,u)$ is holomorphic by $u$ and antiholomorphic by $t$,
where $\mbox{ }^jz={{\mbox{ }^jt\quad \mbox{ }^ju}
\choose {-\mbox{ }^j{\bar u}\quad \mbox{ }^j{\bar t}}}$,
$\mbox{ }^jz\in \bf H$, $z=(\mbox{ }^1z,...,\mbox{ }^nz)\in U$.
Therefore, $N(f)=N(f_{1,1})\cap N(f_{1,2})$, consequently,
$U\setminus N(f)=(U\setminus N(f_{1,1}))\cup (U\setminus N(f_{1,2})).$
Then from Corollary 1.2.4 \cite{henlei} for complex holomorphic functions
$(i)$ follows.
\par $(ii)$. Suppose that $N(f)$ is compact.
In view of $(i)$ and Theorem 2.12
the function $1/f$ can be quaternion holomorphically extended
on $N(f)$. This is a contardiction, since $f=0$ on $N(f)$.
\par {\bf 2.14.1. Note.} Corollary $2.14$ is not true for arbitrary
quaternion holomorphic functions, for example,
$f(\mbox{ }^1z,\mbox{ }^2z)=f_1(\mbox{ }^1z)f_2(\mbox{ }^2z)$
on $B({\bf H^2},0,2)$, where $f_1(\mbox{ }^1z):=
\mbox{ }^1z\mbox{ }^1{\tilde z}-r_1$, $f_2(\mbox{ }^2z):=
\mbox{ }^2z\mbox{ }^2{\tilde z}-r_2$, $0<r_1$, $0<r_2$, $r_1^2+r_2^2<4$.
\par {\bf 2.15. Theorem.} {\it Let $U$ be an open subset
in $\bf H^n$, $f_1$,...,$f_n$ be infinite differentiable
(by real variables) functions on $U$ and Conditions
$2.11.(i)$ are satisfied in $U$. Then for each open
bounded polytor $P=P_1\times ...\times P_n$ such that
$cl (P)$ is a subset in $U$,
there exists a function $u$ infinite differentiable (by real variables)
on $P$ and satisfying Conditions $2.11.(ii)$ on $P$.}
\par {\bf Proof.} Suppose that the theorem is true for
$f_{m+1}=...=f_n=0$ on $U$. The case $m=0$ is trivial.
Assume that the theorem is proved for $m-1$.
Consider $U'={U'}_1\times ...\times {U'}_n$ and
$U"={U"}_1\times ...\times {U"}_n$ open polytors in $\bf H^n$
such that $P\subset cl (P)\subset {U"}\subset cl ({U"})
\subset {U'}\subset cl (U')\subset U$. Take an infinite
differentiable (by real variables) function $\chi $ on
${U'}_m$ with compact support such that $\chi |_{{U"}_m}=1$,
$\chi =0$ in a neighbourhood of ${\bf H}\setminus {U'}_m$.
There exists a function
$$\eta (z):=-(2\pi )^{-3}\int_{{U'}_m}\{ [ \chi (\zeta )
({\hat f}_m(\mbox{ }^1z,...,\mbox{ }^{m-1}z,\zeta _1,
\mbox{ }^{m+1}z,...,\mbox{ }^nz).d{\tilde \zeta }_4) \wedge \nu ,$$
where a differential form $\nu $ is given in \S 2.3
with $\zeta _1, \zeta _2, \zeta _3 \in {U'}_m$ and
$\mbox{ }^mz$ here for $\nu $ instead of $z$ in \S 2.3.
By changing of variables as in \S 2.3 we get
$$\eta (z):=-(2\pi )^{-3}\int_{\bf H}\chi (\zeta _1+z)
({\hat f}_m(\mbox{ }^1z,...,\mbox{ }^{m-1}z,
\zeta _1+\mbox{  }^mz,\mbox{ }^{m+1}z,...,\mbox{ }^nz
).d{\tilde \zeta }_4)\wedge \psi ,\mbox{ where}$$
$$\psi :=(\partial _{\zeta _1} Ln (\zeta _1-\zeta _2))M_1^{-1}
(\partial _{\zeta _2} Ln (\zeta _2-\zeta _3))M_2^{-1}
(\partial _{\zeta _3} Ln (\zeta _3)M_3^{-1}.$$
Consequently, $\partial \eta /\partial {\tilde z}_m={\hat f}_m$ in
${U"}$. The final part of the proof is analogous to that
of Theorem 1.2.5 \cite{henlei}.
\par {\bf 2.16. Definition.} Let $W$ be an open subset in $\bf H^n$
and for each open subsets $U$ and $V$ in $\bf H^n$ such that
\par $(i)$ $\emptyset \ne U\subset V\cap W\ne V$ and
\par $(ii)$ $V$ is connected  \\
there exists a quaternion holomorphic (right superlinear
superdifferentiable, in short RSS, correspondingly) function
$f$ in $W$ such that there does not exist any quaternion
holomorphic (RSS) function $g$ in $V$ such that $g=f$ in $U$.
Then $W$ is called a domain of quaternion (RSS, respectively)
holomorphy. Sets of quaternion holomorphic (RSS) functions in $W$
are denoted by ${\cal H}(W)$ (${\cal H}_{RSL}(W)$ respectively).
\par {\bf 2.17. Definition.} Suppose that $W$ is an open subset
in $\bf H^n$ and $K$ is a compact subset of $W$, then
\par $(i)$ ${\hat K}^{\cal H}_W:= \{ z\in W: |f(z)| \le
\sup_{\zeta \in K}\| {\hat f}(\zeta ) \| $
$\mbox{for each}$ $f\in {\cal H}(W) \} $;
\par $(ii)$ ${\hat K}^{{\cal H}_{RSL}}_W:= \{ z\in W: |f(z)| \le
\sup_{\zeta \in K}|f(\zeta )|$ $\mbox{for each}$
$f\in {\cal H}_{RSL}(W) \} $;  \\
these sets are called the ${\cal H}(W)$-convex hull of $K$
and the ${\cal H}_{RSL}(W)$-convex hull of $K$ respectively,
where $\| {\hat f}(\zeta )\| :=\sup_{h\in {\bf H^n}, |h|\le 1}
|{\hat f}(\zeta ).h|$.
If $K={\hat K}^{\cal H}_W$ or $K={\hat K}^{{\cal H}_{RSL}}_W$,
then $K$ is called ${\cal H}(W)$-convex or ${\cal H}_{RSL}(W)$-convex
correspondingly.
\par {\bf 2.18. Proposition.} {\it For each compact set $K$
in $\bf H^n$, the ${\cal H}({\bf H^n})$-hull and ${\cal H}_{RSL}(
{\bf H^n})$-hull of $K$ are contained in the $\bf R$-convex hull
of $K$.}
\par {\bf Proof.} $I.$ Consider at first the ${\cal H}({\bf H^n})$-hull
of $K$. Each $z\in \bf H^n$ can be written in the form
$z=(\mbox{ }^1z,...,\mbox{ }^nz)$, $\mbox{ }^jz\in \bf H$,
$\mbox{ }^jz=\sum_{l=1}^4x_{l,j}S_l$,
where $x_{l,j}=x_{l,j}(z)\in \bf R$, $S_l\in \{ I,J,K,L \} $.
If $w\in \bf H^n$, $w\notin co_{\bf R}(K)$, then
there are $y_1,...,y_{4n}\in \bf R$ such that $\sum_{j=1}^n\sum_{l=1}^4
x_{l,j}(w)y_{4(j-1)+l}=0$, but \\
$\sum_{j=1}^n\sum_{l=1}^4x_{l,j}(w) y_{4(j-1)+l}<0$ if $z\in K$, where \\
$co_{\bf R}(K):=\{ z\in {\bf H^n}:$ $\mbox{there are}$ $a_1,...,a_s\in \bf R$
$\mbox{and}$ $v_1,...,v_s\in K$ $\mbox{such that}$
$z=a_1v_1+...+a_kv_k \} $ denotes a $\bf R$-convex
hull of $K$ in $\bf H^n$.
Put $\zeta _j=\sum_{l,j}y_{4(j-1)+l}S_l$, then
$f(z):=\exp (\sum_{j=1}^nz_j{\tilde \zeta }_j)$
is the quaternion holomorphic function in $\bf H^n$ such that
$|f(z)|<1$ for each $z\in K$ and $|f(w)|=1$ for the marked
point $w$ above (see Corollary 3.3 \cite{luoyst}), since
$J^2=K^2=L^2=-I$. From $\| {\hat f}(\zeta ) \| \ge
|f(\zeta )|$ the first statement follows.
\par $II.$ Consider now the ${\cal H}_{RSL}({\bf H^n})$-hull
of $K$. Each $f\in {\cal H}_{RSL}(W)$ has the form
$f={ {f_{1,1}\quad f_{1,2}}\choose {-{\bar f}_{1,2}\quad {\bar f}_{2,2}} }$
such that $f_{1,1}$ and $f_{1,2}$ are complex holomorphic
by complex variables $\mbox{ }^jt$ and $\mbox{ }^ju$
respectively and antiholomorphic
by complex variables $\mbox{ }^ju$ and $\mbox{ }^jt$ correspondingly
(see Proposition 2.2 \cite{luoyst}). The set $K$
has projection $K_{1,1}$ and $K_{1,2}$ on the
complex subspaces $\bf C^n$ corresponding to variables
$\mbox{ }^1t,...,\mbox{ }^nt$ and $\mbox{ }^1u,...,\mbox{ }^nu$ respectively.
Therefore, $({\hat K}^{{\cal H}_{RSL}}_{\bf H^n})_{1,l}
\subset {{\hat K}_{1,l,\bf C^n}}^{\cal O}$
for $l=1$ and $l=2$, where ${{\hat K}_{1,l\bf C^n}}^{\cal O}$
denotes the complex holomorphic hull of $K_{1,l}$ in $\bf C^n$.
In view of Proposition 1.3.3
${{\hat K}_{1,l,\bf C^n}}^{\cal O}\subset co_{\bf R}(K_{1,l})$,
hence ${\hat K}^{{\cal H}_{RSL}}_{\bf H^n}\subset co_{\bf R}(K)$.
\par {\bf 2.18.1. Note.} Due to Proposition $2.18$ above
Corollary $1.3.4$ \cite{henlei} can be transferred on $\cal H$ and
${\cal H}_{RSL}$ for $\bf H^n$ instead of $\bf C^n$.
Quaternion versions of Theorems $1.3.5,7,11$, Corollaries
$1.3.6,8,9,10,13$ and Definition $1.3.12$ are true in the
${\cal H}_{RSL}$-class of functions instead of complex holomorphic
functions.
\section{Integral representations of functions of quaternion
variables}
{\bf 3.1. Definitions and Notations.} Consider an $\bf H$-valued
function on $\bf H^n$ such that \\
$(i)\quad (\zeta ,\zeta )=ae$ with
$a\ge 0$ and $(\zeta ,\zeta )=0$ if and only if $\zeta =0$, \\
$(ii)\quad (\zeta ,z+\xi )=(\zeta ,z)+(\zeta ,\xi )$, \\
$(iii)\quad (\zeta +\xi ,z)=(\zeta ,z)+(\xi ,z)$, \\
$(iv)\quad (\alpha \zeta ,z\beta )={\tilde \alpha }(\zeta ,z)\beta $, \\
$(v)\quad (\zeta ,z)^{\tilde .}=(z,\zeta )$ for each $\zeta , \xi $ and
$z\in \bf H^n$, $\alpha $ and $\beta \in \bf H$, $n\in \bf N$.
Then this function is called the scalar product in $\bf H^n$.
The corresponding norm is: \\
$(vi)\quad |\zeta |=\{ (\zeta ,\zeta ) \} ^{1/2}$.
In particular, it is possible to take the canonical scalar product: \\
$(vii)\quad <\zeta ;z>:=
(\zeta ,z)=\sum_{l=1}^n\mbox{ }^l{\tilde \zeta }\mbox{ }^lz$,
where $z=(\mbox{ }^1z,...,\mbox{ }^nz)$, $\mbox{ }^lz\in \bf H$.
\par Consider differential forms on $\bf H$: \\
$(1)\quad \omega _1(\zeta )={\tilde \zeta }d{\tilde \zeta },$ 
$\omega _1(\zeta -z)=({\tilde \zeta }-{\tilde z})d{\tilde \zeta },$
$\omega _1(\zeta ,z)=({\tilde \zeta }-{\tilde z})(d{\tilde \zeta }
-d{\tilde z}),$ \\
$(2)\quad \nu _1(\zeta )=(d{\tilde \zeta }) ({\tilde \zeta }),$ 
$\nu _1(\zeta -z)=(d{\tilde \zeta })({\tilde \zeta }-{\tilde z}),$
$\nu _1(\zeta ,z)=(d{\tilde \zeta }-d{\tilde z})({\tilde \zeta }-
{\tilde z}),$ \\
$(3)\quad \omega _2(\zeta )=(jd{\tilde \zeta }j)\wedge (jd\zeta j),$
$\omega _2(\zeta ,z)=j(d{\tilde \zeta }-d{\tilde z})j\wedge j
(d\zeta -dz)j,$ \\
$(4)\quad \nu _2(\zeta )=d{\tilde \zeta }\wedge d{\tilde \zeta },$
$\nu _2(\zeta ,z)=(d{\tilde \zeta }-d{\tilde z})\wedge
(d{\tilde \zeta }-d{\tilde z}),$ \\
$(5)\quad \omega _4(\zeta )=\nu _2(\zeta )\wedge \omega _2(\zeta ),$ \\
$(6)\quad \omega _4(\zeta ,z)=\nu _2(\zeta ,z)\wedge \omega _2(\zeta ),$ \\
$(7)\quad {\bar \omega }_4(\zeta ,z)=\nu _2(\zeta ,z)\wedge
\omega _2(\zeta ,z).$ \\
With the help of them construct differential forms on $\bf H^n$: \\
$$(8)\quad \theta _z(\zeta ):=(2n-1)! (2\pi )^{-2n}
|\zeta -z|^{-4n} \sum_{s=1}^n \{ \omega _4(\mbox{ }^1\zeta )\wedge ...$$
$$\wedge \omega _4(\mbox{ }^{s-1}\zeta )\wedge 
\omega _1(\mbox{ }^s\zeta -\mbox{ }^sz)\wedge \omega _2(\mbox{ }^s\zeta )
\wedge \omega _4(\mbox{ }^{s+1}\zeta )\wedge ... \wedge
\omega _4(\mbox{ }^n\zeta ) $$
$$+ \omega _4(\mbox{ }^1\zeta )\wedge ...
\wedge \omega _4(\mbox{ }^{s-1}\zeta )\wedge 
\nu _1(\mbox{ }^s\zeta -\mbox{ }^sz)\wedge \omega _2(\mbox{ }^s\zeta )
\wedge \omega _4(\mbox{ }^{s+1}\zeta )\wedge
... \wedge \omega _4(\mbox{ }^n\zeta ) \} ;$$
$$(9)\quad \theta (\zeta ,z):=(2n-1)! (2\pi )^{-2n} |\zeta -z|^{-4n}
\sum_{s=1}^n \{ \omega _4(\mbox{ }^1\zeta ,\mbox{ }^1z) \wedge ...  $$
$$\wedge \omega _4(\mbox{ }^{s-1}\zeta ,\mbox{ }^{s-1}z)\wedge 
\omega _1(\mbox{ }^s\zeta ,\mbox{ }^sz)\wedge \omega _2(\mbox{ }^s\zeta )
\wedge \omega _4(\mbox{ }^{s+1}\zeta ,\mbox{ }^{s+1}z)\wedge ...
\wedge \omega _4(\mbox{ }^n\zeta ,\mbox{ }^nz)
+ \omega _4(\mbox{ }^1\zeta ,\mbox{ }^1z)$$
$$\wedge ...\wedge \omega _4(\mbox{ }^{s-1}\zeta ,\mbox{ }^{s-1}z)\wedge 
\nu _1(\mbox{ }^s\zeta ,\mbox{ }^sz)\wedge \omega _2(\mbox{ }^s\zeta )
\wedge \omega _4(\mbox{ }^{s+1}\zeta ,\mbox{ }^{s+1}z)\wedge
... \wedge \omega _4(\mbox{ }^n\zeta ,\mbox{ }^nz) \} ;$$
$$(10)\quad {\bar \theta }(\zeta ,z):=(2n-1)! (2\pi )^{-2n}
|\zeta -z|^{-4n} \sum_{s=1}^n
\{ {\bar \omega }_4(\mbox{ }^1\zeta ,\mbox{ }^1z)\wedge ... $$
$$\wedge {\bar \omega }_4(\mbox{ }^{s-1}\zeta ,\mbox{ }^{s-1}z)\wedge 
\omega _1(\mbox{ }^s\zeta ,\mbox{ }^sz)\wedge \omega _2(\mbox{ }^s\zeta ,
\mbox{ }^sz) \wedge {\bar \omega }_4(\mbox{ }^{s+1}\zeta ,\mbox{ }^{s+1}z)
\wedge ... \wedge {\bar \omega }_4(\mbox{ }^n\zeta ,\mbox{ }^nz) +$$
$$ {\bar \omega }_4(\mbox{ }^1\zeta ,\mbox{ }^1z)\wedge ...
\wedge {\bar \omega }_4(\mbox{ }^{s-1}\zeta ,\mbox{ }^{s-1}z)\wedge 
\nu _1(\mbox{ }^s\zeta ,\mbox{ }^sz)\wedge \omega _2(\mbox{ }^s\zeta ,
\mbox{ }^sz) \wedge {\bar \omega }_4(\mbox{ }^{s+1}\zeta ,
\mbox{ }^{s+1}z)\wedge ... \wedge {\bar \omega }_4(\mbox{ }^n\zeta ,
\mbox{ }^nz) \} ,$$
where $\zeta $ and $z\in \bf H^n$.
If $U$ is an open subset in $\bf H^n$ and $f$ is a bounded
quaternion differential form on $U$, then by the definition:
$$(11)\quad ({\cal B}_Uf)(z):=\int_{\zeta \in U}f(\zeta )
\wedge \theta (\zeta ,z)$$
for each $z\in \bf H^n$. If in addition $U$ is with a piecewise
$C^1$-boundary (by the corresponding real variables)
and $f$ is a bounded differential form on $\partial U$,
then by the definition:
$$(12)\quad ({\cal B}_{\partial U}f)(z):=
\int_{\zeta \in \partial U}f(\zeta )\wedge \theta (\zeta ,z)$$
for each $z\in \bf H^n$.
\par {\bf 3.2. Theorem.} {\it Let  $U$ be an open subset in $\bf H^n$
with piecewise $C^1$-boundary $\partial U$. Suppose that $f$ is a
continuous function on $cl (U)$ and ${\tilde \partial }f$
is continuous on $U$ in the sense of distributions and has a
continuous extension on $cl (U)$. Then
$$(1)\quad f={\cal B}_{\partial U}f-{\cal B}_U{\tilde \partial f}
\mbox{ on } U,$$ where ${\cal B}_U$ and ${\cal B}_{\partial U}$
are the quaternion integral operators given by Equations $3.1.(11,12)$.}
\par {\bf Proof.} The differential form $\theta (\zeta ,z)$
has the decomposition
$$(2)\quad \theta (\zeta ,z)=\sum_{q=0}^{2n-1}\Upsilon _q(\zeta ,z),$$
where $\Upsilon _q(\zeta ,z)$ is the quaternion differential form
with all terms of degree $4n-q-1$ by $\zeta $ and $\tilde \zeta $
and their multiples on quaternion constants and of degree $q$ by $z$
and $\tilde z$ and their multiples on quaternion constants.
The differential form $f(\zeta )$ has the decomposition
$$(3)\quad f(\zeta )=\sum_{r=0}^mf_r(\zeta ),$$ where
$m=deg (f)$, $f_r(\zeta )$ is with all terms of degree $r$ by $\zeta $
and ${\tilde \zeta }$ and their multiples on quaternion constants.
Then $f_r\wedge \Upsilon _q=0$, when $r>q+1.$
By the definition of integration $\int_{\zeta \in U}f_r(\zeta )
\wedge \Upsilon _q(\zeta ,z)=0$ for $r<q+1.$
If $f$ is a function, then $\int_{\zeta \in \partial U}
f(\zeta )\Upsilon _q(\zeta ,z)=0$ for each $q>0$, since $\partial U$
has the dimension $4n-1$, hence
$$(4)\quad ({\cal B}_{\partial U}f)(z)=\int_{\zeta \in \partial U}
f(\zeta )\theta _z(\zeta ),$$
since $\Upsilon _0(\zeta ,z)=\theta _z(\zeta )$.
If $f$ is a $1$-form, then  $\int_{\zeta \in U}f(\zeta ) \wedge
\Upsilon _q(\zeta ,z)=0$ for each $q>0$, since $U$ has the dimension $4n$,
consequently,
$$(5)\quad ({\cal B}_Uf)(z)=\int_{\zeta \in U}
f(\zeta )\wedge \theta _z(\zeta ).$$
Write $\xi \in \bf H$ in the form $\xi =\alpha e+\beta j$,
then ${\tilde \xi }={\bar \alpha }-\beta j$,
where $\alpha =\alpha _0+\alpha _i{\bf i}$ and $\beta =
\beta _0+\beta _i{\bf i}\in \bf C$, $\alpha _0, \alpha _i,
\beta _0$ and $\beta _i\in \bf R$, ${\bf i}:=(-1)^{1/2}$.
There is the identity $\beta j=j{\bar \beta }$, since $ij=-ji=k$,
where ${\bar \beta }=\beta _0-\beta _i{\bf i}$,
${\bf H}={\bf R}e\oplus {\bf R}i\oplus {\bf R}j\oplus {\bf R}k$,
$i^2=j^2=k^2=-e$, $i={\bf i}e$. Then \\
$(6)\quad d\beta j\wedge d{\bar \beta }j=0,$
$\nu _2(\xi )=-d{\bar \alpha }\wedge d\beta j-
d\beta j\wedge d{\bar \alpha }-d\beta \wedge d{\bar \beta }e$, \\
$(7)\quad \omega _2(\xi )=d\alpha \wedge d{\bar \alpha }e+
2d\alpha \wedge d{\bar \beta }j - d\beta \wedge d{\bar \beta }e,$ \\
$(8)\quad d\xi \wedge d{\tilde \xi }\wedge \omega _2(\xi )=0$,
since $d\xi \wedge d{\tilde \xi }=d\alpha \wedge d{\bar \alpha }e-
2d\alpha \wedge d\beta j + d\beta \wedge d{\bar \beta }e$, \\
$(9)\quad \nu _2(\xi )\wedge \omega _2(\xi )=-d\alpha \wedge d{\bar \alpha }
\wedge d\beta \wedge d{\bar \beta }e=4d\alpha _0\wedge d\alpha _i\wedge
d\beta _0\wedge d\beta _ie$. Therefore, $\nu _2(\xi )\wedge
\omega _2(\xi )/4$ plays the role of the volume element in $\bf H$. Hence \\
$(10)\quad d_{\zeta }(|\zeta -z|^{4n}\theta _z(\zeta ))=
[(2n-1)!(2\pi )^{-2n}](2n)4^nd\mbox{ }^1\alpha _0\wedge d\mbox{ }^1\alpha _i
\wedge d\mbox{ }^1\beta _0\wedge d\mbox{ }^1\beta _i\wedge ...
\wedge d\mbox{ }^n\alpha _0\wedge d\mbox{ }^n\alpha _i
\wedge d\mbox{ }^n\beta _0\wedge d\mbox{ }^n\beta _i,$
since $d_{\zeta }=\partial _{\zeta }+\partial _{\tilde z}$
(see Formula $(2.15)$ \cite{luoyst}).
In view of Proposition $1.7.1$ \cite{henlei} and Formulas
$(6-9)$, $3.1.(vi,vii)$ above the differential form $\theta _z(\zeta )$
is closed in $U\setminus \{ z \} $.
\par There exists $\epsilon _0>0$ such that for each
$0<\epsilon <\epsilon _0$ the ball $B({\bf H^n},z,\epsilon ):=
\{ \zeta \in {\bf H^n}: |\zeta -z|\le \epsilon \} $ and hence
the sphere $S({\bf H^n},z,\epsilon ):=
\{ \zeta \in {\bf H^n}: |\zeta -z|=\epsilon \} =\partial B({\bf H^n},
z,\epsilon )$ are contained in $U$. Apply the Stoke's formula
for matrix-valued functions and differential forms componentwise,
then \\
$(11)\quad \int_{S({\bf H^n},z,\epsilon )}f(\zeta )\theta _z(\zeta )=
\int_{\partial U}f(\zeta )\theta _z(\zeta )-\int_{U_{\epsilon }}
[df(\zeta )]\wedge \theta _z(\zeta )$,
where $U_{\epsilon }:=U\setminus B({\bf H^n},z,\epsilon )$,
$0<\epsilon <\epsilon _0$.
\par There are identities:
$d\xi \wedge jd\zeta =d\xi j\wedge d\zeta $ and
$(\xi d\zeta )^{\tilde .}=[d{\tilde \zeta }]\tilde \xi $ for each
$\xi , \zeta \in \bf H$. Then from Identity $(8)$ it follows, that \\
$(i)\quad d\zeta \wedge jd{\tilde \zeta }\wedge
d{\tilde \zeta }\wedge d\zeta =0,$ \\
$(ii)\quad d\zeta \wedge d{\tilde \zeta }\wedge
jd{\tilde \zeta }\wedge d\zeta =0,$ \\
$(iii)\quad d\zeta \wedge d{\tilde \zeta }\wedge
d{\tilde \zeta }\wedge jd\zeta =0,$
since $j^2=-e$ and ${\bf R}e$ is the centre of the quaternion
algebra $\bf H$, $\alpha $ and $\beta \in \bf C$ commute with $d\alpha $,
$d{\bar \alpha },$ $d\beta $ and $d{\bar \beta }$, where
$\zeta =\alpha e+\beta j$.
From $(i-iii)$ with the help of automorphisms $\zeta \mapsto
j\zeta $ and $\zeta \mapsto \zeta j$ it follows, that \\
$(iv)\quad d\zeta \wedge d{\tilde \zeta }\wedge
d{\tilde \zeta }\wedge d\zeta =0,$  \\
$(v)\quad d\zeta \wedge jd{\tilde \zeta }\wedge
jd{\tilde \zeta }\wedge d\zeta =0,$ \\
$(vi)\quad d\zeta \wedge jd{\tilde \zeta } \wedge
jd{\tilde \zeta }j\wedge jd\zeta j=0.$
Therefore, from $(8)$ and $(vi)$ it follows, that
\par $(12)\quad {\cal B}_Udf={\cal B}_U{\tilde \partial }f$,
since $df=\partial f+{\tilde \partial f}$, where $\partial f(\zeta )=
(\partial f(\zeta )/\partial \zeta ).d\zeta $,
${\tilde \partial }f(\zeta )=(\partial f(\zeta )/\partial {\tilde \zeta }).
d{\tilde \zeta }$, $f(\zeta )=f(\zeta ,{\tilde \zeta })$
is the abbreviated notation.
\par In view of Formula $(10)$ and the Stoke's formula: \\
$(13)\quad \int_{S({\bf H^n},z,\epsilon )}\theta _z(\zeta )=
(2n-1)! (2\pi )^{-2n}[4^n\epsilon ^{-4n}2n]\int_{B({\bf H^n},z,\epsilon )}
(dV)e=e$, where $dV$ is the standard volume element of the Euclidean space
$\bf R^{4n}$. From Formula $(13)$ it follows, that \\
$\lim_{\epsilon \to 0}\int_{S({\bf H^n},z,\epsilon )}
f(\zeta )\theta _z(\zeta )=f(z)$, since  \\
$\int_{S({\bf H^n},z,\epsilon )}(f(\zeta )-f(z))\theta _z(\zeta )=
(2n-1)!(2\pi )^{-2n}\epsilon ^{-4n+1}\int_{S({\bf H^n},z,\epsilon )}
(f(\zeta )-f(z))[|\zeta -z|^{4n-1}\theta _z(\zeta )]$.
The form $[|\zeta -z|^{4n-1}\theta _z(\zeta )]$ is bounded on $U$,
consequently, $|\int_{S({\bf H^n},z,\epsilon )}
(f(\zeta )-f(z))\theta _z(\zeta )|\le C_1\max \{ |f(\zeta )-f(z)|:
\zeta \in B({\bf H^n},z,\epsilon ) \} $, where $C_1$ is a positive
constant independent of $f$ and $\epsilon $ for each
$0<\epsilon <\epsilon _0$.
Therefore, Formula $(1)$ follows from Formula $(11)$
by taking the limit when $\epsilon >0$ tends to zero
and using Identity $(12)$.
\par {\bf 3.3. Corollary.} {\it Let $U$ be an open set in $\bf H^n$
and $f$ be a continuous function on $cl (U)$ and quaternion
holomorphic on $U$. Then  \\
$(1)\quad f={\cal B}_{\partial U}f$ on $U$, \\
where ${\cal B}_U$ and ${\cal B}_{\partial U}$ are
the integral operators given by Equations $3.1.(11,12)$.}
\par {\bf Proof.} From ${\tilde \partial f}=0$, since
$\partial f(\zeta )/\partial {\tilde \zeta }=0$,
and Formula $3.2.(1)$ it follows Formula $3.3.(1)$.
\par {\bf 3.4. Definitions and Notations.}
Suppose that $U$ is a bounded open subset in $\bf H^n$
and $\psi (\zeta ,z)$ be a quaternion-valued $C^1$-function
(by the corresponding real variables) defined on $V\times U$,
where $V$ is a neighbourhood of $\partial U$ in $\bf H^n$, such that \\
$(1)\quad <\psi (\zeta ,z);\zeta -z>\ne 0$ for each
$(\zeta ,z)\in {\partial U}\times U$. Then $\psi $ is called
a quaternion boundary distinguishing map.
Consider the function: \\
$(2)\quad \eta ^{\psi }(\zeta ,z,\lambda ):=
\lambda (\zeta -z) <\zeta -z;\zeta -z>^{-1}$ \\
$+(1-\lambda )\psi (\zeta ,z)<\zeta -z;\psi (\zeta ,z)>^{-1}$, \\
(see Formula $3.1.(vii)$) and the differential forms: \\
$(3)\quad \omega _1(\mbox{ }^s{\tilde \eta }^{\psi }(\zeta ,z,\lambda )):=
\mbox{ }^s{\tilde \eta }^{\psi }(\zeta ,z,\lambda )(
{\tilde \partial }_{\mbox{ }^s\zeta ,\mbox{ }^sz}
+d\lambda )\mbox{ }^s{\tilde \eta }^{\psi }(\zeta ,z,\lambda )$, \\
$\omega _1(\mbox{ }^s{\tilde \eta }^{\psi }(\zeta ,z,0)):=
\mbox{ }^s{\tilde \eta }^{\psi }(\zeta ,z,0)
{\tilde \partial }_{\mbox{ }^s\zeta ,\mbox{ }^sz}
\mbox{ }^s{\tilde \eta }^{\psi }(\zeta ,z,0)$, \\
$(4)\quad \nu _1(\mbox{ }^s{\tilde \eta }^{\psi }(\zeta ,z,\lambda )):=
[({\tilde \partial }_{\mbox{ }^s\zeta ,\mbox{ }^sz}
+d\lambda )\mbox{ }^s{\tilde \eta }^{\psi }(\zeta ,z,\lambda )]
\mbox{ }^s{\tilde \eta }^{\psi }(\zeta ,z,\lambda )$, \\
$\nu _1(\mbox{ }^s{\tilde \eta }^{\psi }(\zeta ,z,0)):=
[{\tilde \partial }_{\mbox{ }^s\zeta ,\mbox{ }^sz}
\mbox{ }^s{\tilde \eta }^{\psi }(\zeta ,z,0)]
\mbox{ }^s{\tilde \eta }^{\psi }(\zeta ,z,0)$, \\
$(5)\quad \nu _2(\mbox{ }^s{\tilde \eta }^{\psi }(\zeta ,z,\lambda )):=
[({\tilde \partial }_{\mbox{ }^s\zeta ,\mbox{ }^sz} +d\lambda )
\mbox{ }^s{\tilde \eta }^{\psi }(\zeta ,z,\lambda )] \wedge
[({\tilde \partial }_{\mbox{ }^s\zeta ,\mbox{ }^sz} +d\lambda )
\mbox{ }^s{\tilde \eta }^{\psi }(\zeta ,z,\lambda )]$,
$\nu _2(\mbox{ }^s{\tilde \eta }^{\psi }(\zeta ,z,0)):=
[{\tilde \partial }_{\mbox{ }^s\zeta ,\mbox{ }^sz}
\mbox{ }^s{\tilde \eta }^{\psi }(\zeta ,z,0)] \wedge
[{\tilde \partial }_{\mbox{ }^s\zeta ,\mbox{ }^sz}
\mbox{ }^s{\tilde \eta }^{\psi }(\zeta ,z,0)]$,  \\
analogously to $(3-5)$ are defined $\omega _1(\mbox{ }^s
{\tilde \psi }(\zeta ,z)),$ $\nu _1(\mbox{ }^s{\tilde \psi }
(\zeta ,z)),$ $\nu _2(\mbox{ }^s{\tilde \psi }(\zeta ,z))$
with $\mbox{ }^s{\tilde \psi }(\zeta ,z)$ instead of
$\mbox{ }^s{\tilde \eta }^{\psi }(\zeta ,z,0);$
$$(6)\phi _{\zeta ,z}:=\phi _{\zeta ,z}(\psi (\zeta ,z);\zeta ):=
(2n-1)! (2\pi )^{-2n} <\psi (\zeta ,z); \zeta -z> ^{-2n} $$
$$\sum_{s=1}^n \{ \nu _2( \mbox{ }^1{\tilde \psi }(\zeta ,z))
\wedge \omega _2(\mbox{ }^1\zeta )\wedge ...
\wedge \nu _2(\mbox{ }^{s-1}{\tilde \psi }(\zeta ,z))\wedge
\omega _2(\mbox{ }^{s-1}\zeta )$$
$$\wedge [\omega _1(\mbox{ }^s{\tilde \psi }(\zeta ,z))\wedge
\omega _2(\mbox{ }^s\zeta )+ \nu _1(\mbox{ }^s{\tilde \psi } (\zeta ,z))
\wedge \omega _2(\mbox{ }^s\zeta )] \wedge $$
$$\nu _2(\mbox{ }^{s+1}{\tilde \psi }(\zeta ,z))\wedge
\omega _2(\mbox{ }^{s+1}\zeta ) \wedge  ...\wedge \nu _2(
\mbox{ }^n{\tilde \psi }(\zeta ,z))\wedge \omega _2(\mbox{ }^n\zeta ) \} ;$$
$$(7){\bar \phi }_{\zeta , z, \lambda }:={\bar \phi }_{\zeta , z, \lambda }
(\psi (\zeta ,z); \zeta ):= (2n-1)! (2\pi )^{-2n} $$
$$\sum_{s=1}^n \{ \nu _2(\mbox{ }^1{\tilde \eta }^{\psi }(\zeta ,z,\lambda ))
\wedge \omega _2(\mbox{ }^1\zeta )\wedge ...
\wedge \nu _2(\mbox{ }^{s-1}{\tilde \eta }^{\psi }(\zeta ,z,\lambda ))\wedge
\omega _2(\mbox{ }^{s-1}\zeta )$$
$$\wedge [\omega _1(\mbox{ }^s{\tilde \eta }^{\psi }(\zeta ,z,\lambda ))
\wedge \omega _2(\mbox{ }^s\zeta )+\nu _1 (\mbox{ }^s{\tilde \eta }^{\psi }
(\zeta ,z,\lambda ))\wedge \omega _2(\mbox{ }^s\zeta )] \wedge $$
$$\nu _2(\mbox{ }^{s+1}{\tilde \eta }^{\psi }(\zeta ,z,\lambda ))\wedge
\omega _2(\mbox{ }^{s+1}\zeta )
\wedge ...\wedge \nu _2(\mbox{ }^n{\tilde \eta }^{\psi }
(\zeta ,z,\lambda ))\wedge \omega _2(\mbox{ }^n\zeta ) \} .$$
If $f$ is a bounded differential form on $U$, then define:
$$(8)\quad (L^{\psi }_{\partial U}f)(z):=\int_{\zeta \in \partial U}
f(\zeta ) \wedge \phi _{\zeta ,z}(\psi (\zeta ,z);\zeta ),$$
$$(9)\quad (R^{\psi }_{\partial U}f)(z):=\int_{\zeta \in \partial U,
0\le \lambda \le 1} f(\zeta ) \wedge {\bar \phi }_{\zeta , z, \lambda }(
\psi (\zeta ,z);\zeta ).$$
\par {\bf 3.5. Theorem.} {\it Let $U$ be an open subset in $\bf H^n$
with a piecewise $C^1$-boundary and let $\psi $ be a quaternion boundary
distinguishing map for $U$. Suppose that $f$ is a continuous
mapping $f: cl(U)\to \bf H$ such that ${\tilde \partial }f$
is also continuous on $U$ in the sence of distributions and
has a continuous extension on $cl (U)$. Then
$$(1)\quad f=(L^{\psi }_{\partial U}f)-(R^{\psi }_{\partial U}
{\tilde \partial }f)-(B_U{\tilde \partial }f)\mbox{ on }U,$$
where the quaternion integral operators $B_U$,
$L^{\psi }_{\partial U}$ and $R^{\psi }_{\partial U}$
are given by Equations $3.1.(11)$, $3.4.(8),(9)$.}
\par {\bf Proof.} There is the decomposition:
\par $(2)\quad {\bar \phi }_{\zeta , z, \lambda }=\sum_{q=0}^{2n-1}
\Upsilon ^{\psi }_q (\zeta ,z,\lambda ),$ \\
where $\Upsilon ^{\psi }_q (\zeta ,z,\lambda )$ is a differential
form with all terms of degree $q$ by $z$ and $\tilde z$ and their multiples
on quaternion constants and of degree $(4n-q-1)$
by $(\zeta , \lambda )$ (including $\tilde \zeta $ and multiples
of $\zeta $ and $\tilde \zeta $ on quaternion constants).
A differential form $f$ has Decomposition $3.2.(3)$.
If $\psi (z)$ is a quaternion $z$-superdifferentiable nonzero
function on an open set $V$ in $\bf H^n$, then differentiating
the equality $(\psi (z))(\psi (z))^{-1}=e$ gives
$[d_z(\psi (z))^{-1}].h=-(\psi (z))^{-1}(d_z\psi (z).h)(\psi (z))^{-1}$
for each $z\in V$ and each $h\in \bf H^n$. Then \\
$\int_{\zeta \in \partial U, 0\le \lambda \le 1}
f_r(\zeta )\wedge \Upsilon ^{\psi }_q(\zeta ,z,\lambda )=0$ 
for each $r\ne q+1$, \\
since $dim (\partial U)=4n-1$, $d\lambda \wedge d\lambda =0$ and
$d\lambda $ commutes with each $b\in \bf H$. Therefore,\\
$(3)\quad R^{\psi }_{\partial U}f_r=
\int_{\zeta \in \partial U, 0\le \lambda \le 1}
f_r(\zeta )\wedge \Upsilon ^{\psi }_{r-1}(\zeta ,z,\lambda )
\mbox{ for each } 1\le r \le 2n$
and $R^{\psi }_{\partial U}f_r=0$
for $r=0$ or $r>2n$. In particular, if $f=f_1$, then \\
$(4)\quad R^{\psi }_{\partial U}f_1=
\int_{\zeta \in \partial U, 0\le \lambda \le 1}
f_1(\zeta )\wedge {\bar \phi }_{\zeta , \lambda }
(\psi (\zeta ,z); \zeta ),$ \\
where ${\bar \phi }_{\zeta , \lambda }
(\psi (\zeta ,z); \zeta ),$ is obtained from
${\bar \phi }_{\zeta , z, \lambda } (\psi (\zeta ,z); \zeta )$
by substituting all ${\tilde \partial }_{\mbox{ }^s\zeta ,\mbox{ }^sz}$
in Formulas $3.4.(3-5,7,9)$ on ${\tilde \partial }_{\mbox{ }^s\zeta }$.
On the other hand, with the help of Formulas $3.2.(8), (vi)$
each quaternion external derivative ${\tilde \partial }_{\mbox{ }^s\zeta }$
can be replaced on $d_{\mbox{ }^s\zeta }$ in
${\bar \phi }_{\zeta , \lambda }(\psi (\zeta ,z); \zeta )$
in Formula $(4)$.
For $\phi _{\zeta ,z}$ there is the dezcomposition: \\
$(5)\quad \phi _{\zeta ,z }=\sum_{q=0}^{2n-1}
\Upsilon ^{\psi }_q (\zeta ,z),$
where $\Upsilon ^{\psi }_q(\zeta , z)$ is a differential form
with all terms of degree $q$ by $z$ and $\tilde z$ and their multiples
on quaternion constants and of degree $4n-q-1$ by $\zeta $ and $\tilde z$
and their multiples on quaternion constants. Therefore, \\
$(6)\quad L^{\psi }_{\partial U}f_r=\int_{\zeta \in \partial U}
f_r\wedge \Upsilon ^{\psi }_r(\zeta ,z)$
for each $0\le r\le 2n-1$ \\
and $L^{\psi }_{\partial U}f_r=0$ for $r\ge 2n$.
In particular, for $f=f_0$: \\
$(7)\quad L^{\psi }_{\partial U}f_0= \int_{\zeta \in \partial U} f_0(\zeta )
\phi _{\zeta } (\psi (\zeta ,z); \zeta ),$ \\
where $\phi _{\zeta } (\psi (\zeta ,z)),$ is obtained from
$\phi _{\zeta , z} (\psi (\zeta ,z); \zeta )$
by substituting all ${\tilde \partial }_{\zeta ,z}$
in Formulas $3.4.(3-5,6,8)$ on ${\tilde \partial }_{\zeta }$.
\par In view of Formula $3.2.(1)$ it remains to prove, that
$R^{\psi }_{\partial U}{\tilde \partial }f=
L^{\psi }_{\partial U}f-B_{\partial U}f$ on $U$.
For each $\zeta $ in a neighbourhood of $\partial U$ there
is the identity:
\par $(8)\quad <\eta ^{\psi }(\zeta , z, \lambda ); \zeta -z>=1$
for each $0\le \lambda \le 1$, hence $d_{\zeta ,z,\lambda }
<\eta ^{\psi }(\zeta ,z,\lambda );\zeta -z>=0$.
By Proposition $1.7.1$ \cite{henlei} and Formulas $3.2.(8),(9),(i-vi):$
\par $(9)\quad d_{\zeta ,\lambda }{\bar \phi }_{\zeta , z, \lambda }=0.$
From Identities $3.2.(8), (vi)$ it follows, that
\par $(10)\quad \partial _{\zeta }f\wedge {\bar \phi }_{\zeta , \lambda }=0.$
Therefore, from $(4), (9), (10)$ it follows, that
\par $(11)\quad d_{\zeta ,\lambda }[f(\zeta ){\bar \phi }_{
\zeta , \lambda }]=[{\tilde \partial }_{\zeta }f(\zeta )]\wedge
{\bar \phi }_{\zeta , \lambda },$
since ${\tilde \partial }_{\zeta }f(\zeta )=
\sum_{s=1}^n(\partial  f(\zeta , {\tilde \zeta })/\partial
\mbox{ }^s{\tilde \zeta }).d\mbox{ }^s{\tilde \zeta }.$
In view of Proposition $1.7.3$ \cite{henlei} and Formulas
$3.2.(8),(9),(i-vi);$ $3.4.(1-7)$:
\par $(12)\quad {\bar \phi }_{\zeta , \lambda }|_{\lambda =0}=
\phi _{\zeta }$, ${\bar \phi }_{\zeta , \lambda }|_{\lambda =1}=
\theta _z(\zeta )$. \\
From the Stoke's formula for matrix-valued differential forms,
in particular, for $[f(\zeta ){\bar \phi }_{\zeta ,z,\lambda }
(\psi (\zeta ,z);\zeta )]$ on $\partial U\times [0,1]$
and Formulas $(4), (7), (11), (12)$ above it follows the statement
of this theorem.
\par {\bf 3.6. Corollary.} {\it Let conditions of Theorem 3.5
be satisfied and let $f$ be a quaternion holomorphic function
on $U$, then $f=L^{\psi }_{\partial U}f$ on $U$.}
\par {\bf 3.7. Remark.} For $n=1$ Formula $3.2.(1)$
produces another analog of the Cauchy-Green formula
(see Theorem $2.1$ and Remark $2.1.1$)
without using the quaternion line integrals.
This is caused by the fact that the dimension of $\bf H$
over $\bf R$ is greater, than $2$: $\quad dim_{\bf R}{\bf H}=4$,
that produces new integral relations.
Theorem $3.2$ can be used instead of Theorem $2.1$
to prove theorems $2.3$ and $2.11$ (with differential
forms of Theorem $3.2$ instead of differential forms of Theorem $2.1$).
If $\psi (\zeta ,z)=\zeta -z$, then $L^{\psi }_{\partial U}=
B_{\partial U}$ and $R^{\psi }_{\partial U}=0$, hence
Formula $3.5.(1)$ reduces to Formula $3.2.(1)$.
For a function $f$ or a $1$-form $f$ Formulas $3.2.(4), (5)$
respectively are valid as well for ${\bar \theta }(\zeta ,z)$
instead of $\theta _z(\zeta )$, where $d_{\zeta ,z}{\bar \theta }
(\zeta ,z)=0$ for each $\zeta \ne z$. A choise of $\omega _4$
(see $3.1.(5)$) and the corresponding to it $\omega _1$,
$\nu _1$, $\nu _2$, $\omega _2$ is not unique, for example,
$d{\tilde \zeta }\wedge d\zeta \wedge d\zeta \wedge d\zeta $
may be taken, since it gives up to a multiplier $Ce$,
where $C$ is a real constant, the canonical volume element
in $\bf H$ and $d\zeta \wedge d\zeta \wedge d\zeta \wedge d\zeta =0.$
Formulas $3.2.(1)$ and $3.5.(1)$ for functions of quaternion variables
are the quaternion analogs of the Martinelli-Bochner and the Leray
formulas for functions of complex variables
respectively, where $\psi (\zeta ,z)$ is the quaternion analog of the
Leray complex map (see \S 3.4).
In the quaternion case the algebra of differential forms
bears the additional gradation structure and have another
properties, than in the complex case (see also \S \S 2.8 and 3.7
\cite{luoyst}). Lemma $3.9$ below shows, that quaternion
boundary distinguishing maps exist.
\par {\bf 3.8. Definitions and Notations.}
Let a subset $U$ in $\bf H^n$ be given by:
\par $(1)$ $U:=\{ z\in {\bf H^n}:$ $\rho (z)<0 \} $,
where $\rho $ is a real-valued function such that
there exists a constant $\epsilon _0>0$ for which:
\par $(2)$ $\sum_{l,m=1}^{4n}(\partial ^2\rho (z)/
\partial x_l\partial x_m)t_lt_m\ge \epsilon _0|t|^2$
for each $t\in \bf R^{4n}$, where $z=(\mbox{ }^1z,...,
\mbox{ }^nz)$, $\mbox{ }^lz\in \bf H$, $\mbox{ }^lz=
\sum_{m=1}^4x_{4(l-1)+m}S_m$, $S_1:=e$, $S_2:=i$, $S_3:=j$,
$S_4:=k$, $x_l\in \bf R$. Then $U$ is called a strictly
convex open set (with $C^2$-boundary). Let
\par $(3)$ $w_{\rho } (z):=(\partial \rho (z)/ \partial \mbox{ }^1z,...,
\partial \rho (z)/ \partial \mbox{ }^nz)$ and
$v_{\rho }(z):=\sum_{m=1}^4(w_{\rho }.S_m)S_m,$
where as usually $w_{\rho }.S_m=(d_z\rho (z)).S_m$
is the differential of $\rho $.
\par {\bf 3.9. Lemma.} {\it Let the function $v_{\rho }$ be as in \S 3.8.
Then $v_{\rho }$ is the quaternion boundary distinguishing map for $U$.}
\par {\bf Proof.} Since $S_mS_l=(-1)^{\kappa (S_m)+\kappa (S_l)}S_lS_m$,
where $\kappa (S_1)=0$, $\kappa (S_2)=\kappa (S_3)=\kappa (S_4)=1$, then \\
$<v_{\rho }(\zeta );\zeta -z>+<\zeta -z;v_{\rho }(\zeta )>=2
\sum_{l=1}^{4n}(\partial \rho (\zeta )/\partial x_l) x_l(\zeta -z)$, \\
where $x_l=x_l(\zeta )$ and $x_l(\zeta -z)$
are real coordinates corresponding to $\zeta $ and $\zeta -z$.
By the Taylor's theorem:
$\rho (z)e= \rho (\zeta )e-<v_{\rho }(\zeta );\zeta -z>/2-
<\zeta -z; v_{\rho }(\zeta )>/2+\sum_{l,m=1}^{4n} (\partial ^2\rho (\zeta )/
\partial x_l\partial x_m) x_l(\zeta -z) x_m(\zeta -z)e/2
+o(|\zeta -z|^2)e$. Therefore, there exists a neighbourhood
$V$ of $\partial U$ and $\epsilon _1>0$ such that \\
$(1)\quad (<v_{\rho }(\zeta );\zeta -z>+<\zeta -z; v_{\rho }(\zeta )>)_e/2
\ge \rho (\zeta )-\rho (z)+\epsilon _0 |\zeta -z|^2/4$
for each $\zeta \in V$ and $|\zeta -z|\le \epsilon _1$,
where $a=a_ee+a_ii+a_jj+a_kk$ for each $a\in \bf H$, $a_e$,
$a_i$, $a_j$ and $a_k$ are reals.
If $z\in U$, $\zeta \in \partial U$, $|z-\zeta |\le \epsilon _1$,
then by $(1)$:
$(<v_{\rho }(\zeta ); \zeta -z>+<\zeta -z; v_{\rho }(\zeta )>)_e
\ge -\rho (z)>0$. If $|\zeta -z|>\epsilon _1$, put
$z_1:=(1-\epsilon _1|\zeta -z|^{-1})\zeta +
\epsilon _1|\zeta -z|^{-1}z$, then $\zeta -z_1=\epsilon _1|\zeta -z|^{-1}
(\zeta -z)$, consequently,
$(<v_{\rho }(\zeta ); \zeta -z>+<\zeta -z; v_{\rho }(\zeta )>)_e/2
\ge -\rho (z_1)$. Evidently, $U$ is convex and $z_1\in U$.
\par {\bf 3.10. Theorem.} {\it Let $U$ be a strictly
convex open subset in $\bf H^n$ (see $3.8.(1)$) and let
$f$ be a continuous function on
$U$ with continuous ${\tilde \partial }f$ on $U$ in the sence
of distributions having a continuous extension on $cl (U)$
such that $2.11.(i)$ is satisfied. Then there exists a function
$u$ on $U$ which is a solution of the ${\tilde \partial }$-equation
$2.11.(ii)$.}
\par {\bf Proof.} In proofs of Theorems $2.3$ and $2.11$
take in Formula $3.5.(1)$ $\chi f$ instead of $f$,
which is possible due to Lemma $3.9$, choosing $\psi =
v_{\rho }$ and $supp (\chi )$ as a proper subset of $U$.
Then $L^{\psi }_{\partial U}\chi f=0$ and
$R^{\psi }_{\partial U}\chi f=0$, hence $\chi f=-B_U{\tilde
\partial }\chi f$. For each fixed $z\in U$ a subset $\mbox{ }^lU_{\eta }:=
\{ \xi \in {\bf H}: \rho (\mbox{ }^1z,...,
\mbox{ }^{l-1}z,\xi ,\mbox{ }^{l+1}z,...,\mbox{ }^nz) <0 \} $
is strictly convex in $\bf H$ due to $3.8.(1,2)$, where
$\eta :=(\mbox{ }^1z,...,\mbox{ }^{l-1}z,\mbox{ }^{l+1}z,...,\mbox{ }^nz).$
Apply $3.5.(1)$ by a variable $\xi $ in $\mbox{ }^lU_{\eta }$,
in particular, for $l=1$, for which $v_{\rho }$ by the variable
$\xi $ is the quaternion boundary distinguishing map for
$\mbox{ }^1U_{\eta }$. Therefore,
$u(z):=-B_{\mbox{ }^1U_{\eta }}\mbox{ }^1{\hat f}(\xi ,\eta ).
d{\tilde \xi }$ with $z=(\xi ,\eta )$, $\xi \in \mbox{ }^1U_{\eta }$
solves the problem.
\section{Quaternion manifolds}
\par {\bf 4.1. Definitions and Notations.}
Suppose that $M$ is a quaternion manifold and let
$GL(N,{\bf H})$ be the group of all invertible quaternion
$N\times N$ matrices. Then a quaternion holomorphic vector bundle
$Q$ of quaternion dimension $N$ over $M$ is a $C^{\infty }$-vector
bundle $Q$ over $M$ with the characteristic fibre
$\bf H^n$ together with a quaternion holomorphic atlas
of local trivializations: $g_{a,b}: U_a\cap U_b\to GL(N,{\bf H})$,
where $U_a\cap U_b\ne \emptyset $, $ \{ (U_a,h_a): a \in \Upsilon \}
=At (Q)$, $\bigcup_a U_a=M$, $U_a$ is open in $M$, $h_a:
Q|_{U_a}\to U_a\times {\bf H^n}$ is the bundle isomorphism,
$(z,g_{a,b}(z)v)=h_a\circ h_b^{-1}(z,v)$, $z\in U_a\cap U_b$,
$v\in \bf H^n$. Since $M$ is also the real manifold there exists
the tangent bundle $TM$ such that $T_xM$ is isomorphic with
$\bf H^n$ for each $x\in M$, since $TU_a=U_a\times \bf H^n$ for each $a$,
where $dim_{\bf H}M=n$ is the quaternion dimension of $M$.
If $X$ is a Banach space over $\bf H$ (with left and right
distributivity laws relative to multiplications of vectors
in $X$ on scalars from $\bf H$), then denote by $X^*_q$
the space of all additive $\bf R$-homogeneous functionals on $X$
with values in $\bf H$. Clearly $X^*_q$ is the Banach space over $\bf H$.
Then $T^*M$ with fibres $(H^n_q)^*$ denotes the quaternion
cotangent bundle of $M$ and $\Lambda ^rT^*M$ denotes the vector
bundle whose sections are quaternion $r$-forms on $M$, where
$S_bdx_b\wedge S_adx_a=-(-1)^{\kappa (S_a)+\kappa (S_b)}
S_adx_a\wedge S_bdx_b$ for each $S_a, S_b\in \{ e,i,j,k \} $,
$dz=edx_e+idx_i+jdx_j+kdx_k$, $z\in \bf H$, $x_b\in \bf R$.
\par Quaternion holomorphic Cousin data in $Q$ is a family
$ \{ f_{a,b}: a, b \in \Upsilon \} $ of quaternion holomorphic sections
$f_{a,b}: U_a\cap U_b\to Q$ such that $f_{a,b}+f_{b,l}=f_{a,l}$
in $U_a\cap U_b\cap U_l$ for each $a, b, l \in \Upsilon $.
A finding of a family $\{ f_a: a\in \Upsilon \} $ of quaternion
holomorphic sections $f_a: U_a\to Q$ such that $f_{a,b}=f_a-f_b$ in
$U_a\cap U_b$ for each $a, b \in \Upsilon $ will be called
the quaternion Cousin problem.
\par {\bf 4.2. Theorem.} {\it Let $M$ be a quaternion manifold
and $Q$ be a quaternion holomorphic vector bundle on $M$. Then Conditions
$(i,ii)$ are equivalent:
\par $(i)$ each quaternion holomorphic Cousin problem in $M$
has a solution;
\par $(ii)$ for each quaternion holomorphic section $f$ of $Q$
such that ${\tilde \partial }f=0$ on $M$, there exists a
$C^{\infty }$-section $U$ of $Q$ such that $(\partial u/\partial
{\tilde z})={\hat f}$ on $M$.}
\par {\bf Proof.} $(i)\rightarrow (ii)$. In view of Theorems $2.11$
and $3.10$ there exists an (open) covering $ \{ U_a : a \} $
of $M$ and $C^{\infty }$-sections $u_b: U_b\to Q$ such that
$(\partial u_b/\partial {\tilde z})={\hat f}$ in $U_b$.
Then $(u_b-u_l)$ is quaternion holomorphic in $U_b\cap U_l$
and their family forms quaternion holomorphic Cousin data in $Q$. Put
$u:=u_b-h_b$ on $U_b$, where $u_b-u_l=h_b-h_l$, $h_b:
U_b\to Q$ is a quaternion holomorphic section given by $(i)$.
\par $(ii)\rightarrow (i)$.
Take a $C^{\infty }$-partition of unity $\{ \chi _b: b \} $
subordinated to $ \{ U_b: b \} $ and $c_b:=-\sum_a\chi _af_{a,b}$
on $U_b$, then $f_{l,b}=\sum_a\chi _a(f_{l,a}+f_{a,b})=c_l-c_b$
in $U_l\cap U_b$, hence $(\partial c_l/\partial {\tilde z})
=(\partial c_b/\partial {\tilde z})$ in $U_l\cap U_b$.
By $(ii)$ there exists a $C^{\infty }$-section $u: M\to Q$
with $(\partial u/\partial {\tilde z})=(\partial c_b/\partial {\tilde z})$
on $U_b$ and $h_b:=c_b-u$ on $U_b$ gives the solution.
\par {\bf 4.3. Definitions.}
Suppose $U$ is an open subset in $\bf H$, then a $C^2$-function
$\rho : U\to \bf R$ is called subharmonic (strictly subharmonic)
in $U$ if $\sum_{m=1}^4\partial ^2\rho /\partial x_m^2\ge 0$
($\sum_{m=1}^4\partial ^2\rho /\partial x_m^2>0$ correspondingly)
for each $z=x_1e+x_2i+x_3j+x_4k\in U$, where $x_1,...,x_4\in \bf R$.
If $U$ is an open subset in $\bf H^n$, then a $C^2$-function
$\rho : U\to \bf R$ such that the function $\zeta \mapsto \rho
(v+\zeta w)$ is subharmonic (strictly subharmonic) on its domain
for each $v, w\in \bf H^n$ is called plurisubharmonic (strictly
plurisubharmonic correspondingly) function, where $\zeta \in \bf H$.
\par A $C^p$-function $\rho $ on a quaternion manifold $M$
is called a strictly plurisubharmonic exhausting $C^p$-function
for $M$, $2\le p\in \bf N$, if $\rho $ is a strictly plurisubharmonic
$C^p$-function on $M$ and for each $\alpha \in \bf R$ the set
$ \{ z\in M: \rho (z)<\alpha \} $ is relatively compact in $M$.
\par {\bf 4.4. Theorem.} {\it Let $M$ be a quaternion manifold
with strictly plurisubharmonic exhausting function $\rho $
such that $\rho e$ is a $C^{\omega }_{z,\tilde z}$-function and let
$Q$ be a quaternion holomorphic vector bundle on $M$,
$U_{\alpha }:= \{ z\in M: \rho (z)<\alpha \} $ for $\alpha \in \bf R$.
\par $(i).$ Suppose that $d\rho (z)\ne 0$ for each $z\in \partial
U_{\alpha }$ for a marked $\alpha \in \bf R$. Then every continuous
section $f: cl (U_{\alpha })\to Q$ quaternion holomorphic on $U_{\alpha }$
can be approximated uniformly on $cl (U_{\alpha })$ by quaternion
holomorphic sections of $Q$ on $M$.
\par $(ii).$ For each continuous mapping $f: M\to Q$ such that
${\tilde \partial f}=0$ on $M$ there exists a continuous mapping
$u: M\to Q$ such that $\partial u/\partial {\tilde z}={\hat f}$ on $M$.}
\par {\bf Proof.} For a $C^{\omega }_{z,\tilde z}$-function
$\rho : U\to \bf R$ (that is, $\rho $ is locally analytic in variables
$(z,{\tilde z})$, ${\bf R}={\bf R}e\hookrightarrow \bf H$)
there is the identity: \\
$\sum_{l,m,a,b} (\partial ^2\rho /\partial \mbox{ }^lx_a\partial
\mbox{ }^mx_b)t_{4(l-1)+a}t_{4(l-1)+b}=$ \\
$\sum_{m,l}(\partial ^2 \rho (z)/\partial \mbox{ }^lz\partial
\mbox{ }^mz).((\partial \mbox{ }^lz/\mbox{ }^lx_a)t_{4(l-1)+a},
(\partial \mbox{ }^mz/\mbox{ }^lx_b)t_{4(l-1)+b})$ \\
$= \sum_{m,l=1}^n (\partial ^2 \rho (z)/\partial \mbox{ }^lz\partial
\mbox{ }^mz).(\mbox{ }^l\xi ,\mbox{ }^m{\tilde \xi })$,
since $\partial \mbox{ }^lz/\partial \mbox{ }^lx_a=S_a$,
$\partial \mbox{ }^l{\tilde z}/\partial \mbox{ }^lx_a=
(-1)^{\kappa (S_a)}S_a$, where $\mbox{ }^l\xi =\sum_{a=1}^4
t_{4(l-1)+a}S_a$, $S_1=e$, $S_2=i$, $S_3=j$, $S_4=k$,
$\mbox{ }^lz= \sum_{a=1}^4 \mbox{ }^lx_aS_a$, $t_b\in \bf R$,
$\mbox{ }^lx_a \in \bf R$.  Therefore, a $C^{\omega }_{z,\tilde z}$-function
$\rho $ is strictly plurisubharmonic on $U$ if and only if
\par $(1)\quad \sum_{m,l=1}^n (\partial ^2 \rho (z)
/\partial \mbox{ }^lz\partial
\mbox{ }^m{\tilde z}).(\mbox{ }^l\xi ,\mbox{ }^m{\tilde \xi })>0$
for each $z\in U$ and each $0\ne \xi \in \bf H^n$,
where $\xi =(\mbox{ }^1\xi ,...,\mbox{ }^n\xi )$
(see also \S 2 \cite{luoyst}). Consider a proper compact subset
$A$ in $M$ such that $d\rho (z)\ne 0$ for each $z\in A$.
Then for each $\epsilon >0$ there exists a strictly plurisubharmonic
function $\rho _{\epsilon }: M\to \bf R$ such that $\rho e$ is
a $C^{\omega }_{z,\tilde z}$-function on $M$ and $(i-iii)$ are fulfilled:
\par $(i)$ $\rho -\rho _{\epsilon }$ together with its first
and second derivatives is not greater than $\epsilon $ on $M$;
\par $(ii)$ the set $Crit (\rho _{\epsilon }):=\{ z\in M:
d\rho _{\epsilon }(z)=0 \} $ is discrete in $M$;
\par $(iii)$ $\rho _{\epsilon }=\rho $ on $A$ (see also Lemma
$2.1.2.2$ \cite{henlei} in the complex case).
\par The space $C^{\omega }_z(U,{\bf H})$ is dense in
$C^0(U,{\bf H})$ for each open $U$ in $\bf H^n$ (see \S 2.7
and Theorem $3.28$ in \cite{luoyst}). Suppose $\beta \in \bf R$
and $d\rho (z)\ne 0$ for $z\in \partial U_{\beta }$ and
$f: cl(U_{\beta })\to Q$ is a continuous section quaternion
holomorphic on $U_{\beta }$. Therefore, for each
$\beta \le \alpha <\infty $ if $d\rho (z)\ne 0$ for each
$z\in \partial U_{\alpha }$, then $f$ can be uniformly
approximated on $cl(U_{\beta })$ by continuous sections
on $cl(U_{\alpha })$ that are holomorphic on $U_{\alpha }$.
There exists a sequence $\beta <\alpha _1<\alpha _2<...$ such that
$\lim_l \alpha _l=\infty $ and $d\rho (z)\ne 0$ for each
$z\in \partial U_{\alpha _l}$, since $Crit (\rho )$ is discrete.
For each $\epsilon >0$ there exists a continuous section
$f_l: cl(U_{\alpha _l})\to Q$ such that $f_l$ is quaternion
holomorphic on $U_{\alpha _l}$ and $\| f_{l+1}-f_l \|_{C^0(U_{\alpha _l})}<
\epsilon 2^{-l-1}$ for each $l\in \bf N$, where $f_0:=f$.
Therefore, the sequence $ \{ f_l: l\in {\bf N} \} $ converges
to the quaternion holomorphic section $g: M\to Q$ uniformly
on each compact subset $P$ in $M$ and $\| f-g \|_{C^0(U_{\beta })}<
\epsilon $.
\par The second statement $(ii)$ follows from $(i)$ and
Theorems $2.11, 3.10$, since $Crit (\rho )$ is discrete in $M$
and there exists a sequence of continuous $Q$-valued functions
on $cl (U_{\alpha _l})$ such that $\partial u_l/\partial {\tilde z}=
\hat f$ on $U_{\alpha _l}$, $\bigcup_lU_{\alpha _l}=M$
(see also the complex case in \S 2.12.3  \cite{henlei} mentioning,
that Lemma $2.12.4$ there can be reformulated and proved for
a quaternion manifold $M$ on $\bf H^n$ instead of a complex manifold
on $\bf C^n$).
\par {\bf 4.5. Definitions.} Let $M$ be a quaternion manifold
(see \S 2.10). For a compact subset $G$ in $M$ put:
${\hat G}^{\cal H}_M:=\{ z\in M: |f(z)|\le \sup_{\zeta \in G}
|{\hat f}(\zeta )| \forall f\in {\cal H}(M) \} $.
Such ${\hat G}^{\cal H}_M$ is called the ${\cal H}(M)$-hull of
$G$. If $G={\hat G}^{\cal H}_M$, then $G$ is called ${\cal H}(M)$-convex.
A quaternion manifold $M$ is called quaternion holomorphically
convex if for each compact subset $G$ in $M$ the set
${\hat G}^{\cal H}_M$ is compact.
\par A quaternion manifold $M$ with a countable atlas $At(M)$
having dimension $n$ over $\bf H$ and satisfying $(i,ii)$:
\par $(i)$ $M$ is quaternion holomorphically convex;
\par $(ii)$ for each $z\in M$ there are $\mbox{ }^1f,...,\mbox{ }^nf
\in {\cal H}(M)$ and there exists a neighbourhood $U$ of $z$
such that the map $U\ni \zeta \mapsto (\mbox{ }^1f(\zeta ),...,\mbox{ }^nf
(\zeta ))$ is quaternion biholomorphic (see \S 2.6),
then $M$ is called a quaternion Stein manifold.
\par {\bf 4.6. Remark.} If $M_1$ and $M_2$ are two quaternion
Stein manifolds, then $M_1\times M_2$ is a quaternion Stein manifold.
If $N$ is a closed quaternion submanifold of a quaternion
Stein manifold $M$, then $N$ is also a quaternion Stein manifold.
\par {\bf 4.7. Theorem.} {\it Let $M$ be a quaternion Stein manifold.
Then for each ${\cal H}(M)$-convex compact subset $P$ in $M$,
$P\ne M$ and each neighbourhood $V_P$ of $P$ there exists a strictly
plurisubharmonic exhausting $C^{\omega }_{z,{\tilde z}}$-function
$\rho $ on $M$
such that $\rho <0$ on $P$ and $\rho >0$ on $M\setminus V_P$.}
\par The {\bf proof} of this theorem is analogous to that of
Theorem $2.3.14$ \cite{henlei} in the complex case taking
$\rho (z):=-1+\sum_{l=1}^{\infty }\sum_{k=1}^{N(l)}
f^k_l(z)(f^k_l(z))^{\tilde .}$ for each $z\in M$, where
$f^k_l\in C^{\omega }_{z,\tilde z}$, $M=\bigcup_lP_l$,
$P_l\subset Int (P_{l+1})$ for each $l\in \bf N$, each $P_l$ is
${\cal H}(M)$-convex, $\sum_{k=1}^{N(l)}|f^k_l(z)|^2<2^{-l}$
for each $z\in P_l$, $\sum_{k=1}^{N(l)}|f^k_l(z)|^2>l$
for each $z\in P_{l+2}\setminus U_l$, $U_l:=Int (P_{l+1})$,
$rank [(\partial f^k_l/\partial \mbox{ }^mz)^{k=1,...,N(l)}_{m=1,...,n}]=n$
for each $z\in P_l$.
\par {\bf 4.8. Theorem.} {\it A quaternion manifold $M$ is a
quaternion Stein  manifold if and only if there exists a strictly
plurisubharmonic exhausting $C^{\omega }_{z,\tilde z}$-function
$\rho $ on $M$, then $ \{ z\in M: \rho (z)\le \alpha \} $
is ${\cal H}(M)$-convex for each $\alpha \in \bf R$.}
\par {\bf Proof.} The necessity follows from Theorem $4.7$.
To prove sufficiency suppose $\eta =(\mbox{ }^1\eta ,...,
\mbox{ }^n\eta )$ are quaternion holomorphic coordinates
in a neighbourhood $V_{\xi }$ of $\xi \in M$. Consider
\par $(1)\quad u(z):=2\sum_{l,m=1}^n<v_{\rho }(\xi );\eta (z)-\eta (\xi )>$ \\
$+\sum_{l,m=1}^n(\partial ^2\rho (\xi )/\partial \mbox{ }^l\eta
\partial \mbox{ }^m\eta ).[(\mbox{ }^l\eta (z)-\mbox{ }^l\eta (\xi )),
(\mbox{ }^m\eta (z)-\mbox{ }^m\eta (\xi ))]e$, \\
where $v_{\rho }(\xi )$ is given by $3.8.(3)$. Then $u$ is holomorphic
in $V_{\xi }$ and $u(\xi )=0$. By Lemma $3.9$:
\par $(2)\quad (u(z)+{\tilde u}(z))/2=\rho (z)e-\rho (\xi )e-
\sum_{l,m=1}^n (\partial ^2\rho (\xi )/\partial \mbox{ }^l\eta
\partial \mbox{ }^m\eta ).[(\mbox{ }^l\eta (z)- \mbox{ }^l\eta (\xi )),
(\mbox{ }^m\eta (z)- \mbox{ }^m\eta (\xi ))^{\tilde .}]
+o(|\eta (\xi )-\eta (z)|^2)$.
From the strict plurisubharmonicity of $\rho $ it follows, that
there exists $\beta >0$ and $V_{\xi }$ such that
\par $(3)\quad (u(z)+{\tilde u}(z))/2< \rho (z)-\rho (\xi )
-\beta |\eta (z)-\eta (\xi )|^2$ for each $z\in V_{\xi }$.
Then $\exp (u(\xi ))=1$ and $|\exp (u(z))|<1$ for each
$\xi \ne z\in cl(U_{\alpha })\cap V_{\xi }$ (see Corollary
$3.3$ \cite{luoyst}).
\par If $g: {\bf R}\to \bf H$ is a $C^{\infty }$-function
with compact support, then $g(z{\tilde z})=:\chi (z)$ is a
$C^{\infty }$-function on $\bf H^n$ with compact support
such that $\chi $ is quaternion
$(z,{\tilde z})$-superdifferentiable. Therefore, there exists
a neighbourhood $W_{\xi }\subset V_{\xi }$ of $\xi $ and an infinitely
$(z,{\tilde z})$-superdifferentiable function $\chi $ such that
$\chi |_{W_{\xi }}=1$, $supp (\chi )$ is a proper subset of $V_{\xi }$,
consequently, $\lim_{m\to \infty } \| \exp (mu(z))
(\partial \chi (z)/\partial {\tilde z}) \| _{C^0(U_{\alpha })}=0$,
where $(\partial \chi (z)/\partial {\tilde z})=(\partial \chi (z)/
\partial \mbox{ }^1{\tilde z},...,\partial \chi (z)/\partial \mbox{ }^n
{\tilde z}).$ In view of Theorem $3.10$ there exist continuous functions
$v_m$ on $cl (U_{\alpha })$ such that $(\partial v_m/\partial {\tilde z})=
\exp (mu(z))(\partial \chi /\partial {\tilde z})$ in $U_{\alpha }$
and $\lim_{m\to \infty } \| v_m \|_{C^0(U_{\alpha })}=0$.
\par Put $g_m(z):=\exp (mu(z))\chi (z)-v_m(z)+v_m(\xi )$, hence
$g_m$ is continuous on $cl(U_{\alpha })$ and quaternion holomorphic
on $U_{\alpha }$.
Since $supp (\chi )$ is the proper subset in $V_{\xi }$,
then $g_m(\xi )=1$ for each $m\in \bf N$,
$sup_m \| g_m \|_{C^0(U_{\alpha })}<\infty $ and for each compact subset
$P$ in $cl(U_{\alpha })\setminus \{ \xi \} $ there exists $\lim_m
\| g_m \|_{C^0(P)}=0$. In view of Theorem $4.4.(ii)$ there exists
a sequence of functions $f_m\in {\cal H}(M)$ and $C=const<\infty $
such that $(a)$ $f_m(\xi )=1$ for each $m\in \bf N$;
$(b)$ $\| f_m \|_{C^0(U_{\alpha })}\le C$ for each $m\in \bf N$;
$(c)$ $\lim_{m\to \infty } \| g_m \|_{C^0(P)}=0$ for each
compact subset $P\subset cl(U_{\alpha })\setminus \{ \xi \} $.
\par Consider a quaternion holomorphic function $f$ on a neighbourhood
of $\xi $ such that $f(\xi )=0$. Put $\phi _m:=f\exp (mu)\partial
\chi /\partial {\tilde z}$, then $supp (\phi _m)$ is the proper
subset in $V_{\xi }\setminus W_{\xi }$. In view of Inequality $(3)$
there exists $\delta >0$ such that $\lim_m \| \phi _m \|_{C^0(U_{\alpha
+\delta })}=0$. As in \S 4.4 it is possible to assume, that $Crit (\rho )$
is discrete in $M$. Take $0<\epsilon <\delta $ such that $d\rho \ne 0$
on $\partial U_{\alpha +\epsilon }$. In view of Theorem $4.4.(ii)$
there exists a continuous function $v_m$ on $cl(U_{\alpha +\epsilon })$
such that $\partial v_m/\partial {\tilde z}={\hat \phi }_m$ on
$U_{\alpha +\epsilon }$ and $\lim_m \| v_m \|_{C^0(U_{\alpha +\epsilon })}
=0$. Each $v_m$ is quaternion holomorphic on $W_{\xi }$, since
$\phi _m=0$ on $W_{\xi }$, hence $\lim_m\partial v_m(\xi )=0$.
Since $f(\xi )=u(\xi )=0$ and $\chi =1$ on $W_{\xi }$, then
$\partial g_m(\xi )/\partial \xi =\partial f(\xi )/\partial \xi
-\partial v_m(\xi )/\partial \xi $. In view of Theorem $4.4.(i)$ there
exists $f_m\in {\cal H}(M)$ such that $\| f_m-g_m \|_{C^0(U_{\alpha
+\epsilon })}<m^{-1}$ and inevitably $\lim_m \| \partial f_m(\xi )/
\partial \xi -\partial g_m(\xi )/\partial \xi \| =0$.
\par Let $V_{\xi }$ and $W_{\xi }$ be as above, then there exists
$\delta >0$ such that $(u(z)+{\tilde u}(z))/2<-\delta $ for each
$z\in U_{\alpha +\delta }\cap (V_{\xi })\setminus W_{\xi })$.
Therefore, there exists a branch of the quaternion logarithm
$Ln (u)\in {\cal H}(U_{\alpha +\delta }\cap (V_{\xi }\setminus
cl (W_{\xi }))$ (see \S \S 3.7, 3.8 \cite{luoyst}).
From Theorems $4.2, 4.4$ it follows that each quaternion holomorphic
Cousin problem over $U_{\alpha +\delta }$ has a solution.
Hence $Ln (u)=w_1-w_2$ for suitable $w_1\in {\cal H}(V_{\xi }
\cap U_{\alpha +\delta })$ and $w_2\in {\cal H}(U_{\alpha +\delta }
\setminus cl (W_{\xi })).$ Put $f:=u \exp (-w_1)$ in $U_{\alpha
+\delta }\cap V_{\xi }$ and $f:=\exp (-w_2)$ in $U_{\alpha +\delta }
\setminus cl (W_{\xi })$. Then $f\in {\cal H}(U_{\alpha +\delta })$
and $f(\xi )=0$. In view of Inequality $(3)$ $f(z)\ne 0$ for each
$\xi \ne z\in cl(U_{\alpha })$. Verify now that $cl(U_{\alpha })$
is ${\cal H}(M)$-convex. Consider $\xi \in M\setminus cl(U_{\alpha })$.
Due to \S 4.4 there exists a strictly plurisubharmonic exhausting
$C^{\omega }_{z,{\tilde z}}$-function $\psi $ for $M$ such that
$Crit (\psi )$ is discrete and $U_{\alpha }\subset G_{\psi (\xi )}$,
where $G_{\beta }:= \{ z\in M:$ $\psi (z)<\beta \} $ for $\beta \in \bf R$.
Considering shifts $\psi \mapsto \psi +const$ assume $d\psi (z)\ne 0$
for each $z\in \partial G_{\psi (\xi )}$. From the proof above
it follows, that there exists $f\in {\cal H}(M)$ such that $f(\xi )=
1$ and $|f(z)|<1$ for each $z\in cl(U_{\alpha })$.
\par {\bf 4.8.1. Remark.} With the help of Theorem $4.8$
it is possible to spread certain modifications of Theorems $3.2$
and $3.5$ on quaternion Stein manifolds.
\par {\bf 4.9. Theorem.} {\it Let $N$ be a complex manifold,
then there exists a quaternion manifold $M$ and a
complex  holomorphic embedding $\theta : N\hookrightarrow M$.}
\par {\bf Proof.} Suppose $At (N)=\{ (V_a,\psi _a): a \in \Lambda \} $
is any holomorphic atlas of $N$, where $V_a$ is open in $N$,
$\bigcup_aV_a=N$, $\psi _a: V_a\to \psi _a(V_a)\subset \bf C^n$
is a homeomorphism for each $a$, $n=dim_{\bf C}M\in \bf N$,
$\{ V_a: a \in \Lambda \} $ is a locally finite
covering of $N$, $\psi _b\circ \psi _a^{-1}$ is a holomorphic function
on $\psi _a(V_a\cap V_b)$ for each $a, b\in \Lambda $
such that $V_a\cap V_b\ne \emptyset $.  For each complex holomorphic
function $f$ on an open subset $V$ in $\bf C^n$ there exists a quaternion
holomorphic function $F$ on an open subset $U$ in $\bf H^n$
such that $\pi _{1,1}(U)=V$ and $F_{1,1}|_{Ve}=f|_V$,
where $\pi _{1,1}: {\bf R^n}e\oplus {\bf R^n}i\oplus {\bf R^n}j
\oplus {\bf R^n}k\to {\bf R^n}e\oplus {\bf R^n}i=\bf C^n$
is the projection, $F_{1,1}:=\pi _{1,1}\circ F$
(see Proposition $3.13$ \cite{luoyst} and use a locally finite
covering of $V$ by balls). Therefore, for each two charts
$(V_a,\psi _a)$ and $(V_b,\psi _b)$ with $V_{a,b}:=
V_a\cap V_b\ne \emptyset $ there exists $U_{a,b}$ open in
$\bf H^n$ and a quaternion holomorphic function $\Psi _{b,a}$
such that $\Psi _{b,a}|_{\psi _a(V_{a,b})e}=\psi _{b,a}|_{\psi _a(V_{a,b})}$,
where $\psi _{b,a}:=\psi _b\circ \psi _a^{-1}$,
$\pi _{1,1}(U_{a,b})=\psi _a(V_{a,b})$.
Consider $Q:=\bigoplus_aQ_a$, where $Q_a$ is open in $\bf H^n$,
$\pi _{1,1}(Q_a)=\psi _a(V_a)$ for each $a\in \Lambda $.
The equivalence relation $\cal C$ in the topological space
$\bigoplus_a\psi _a(V_a)$ generated by functions $\psi _{b,a}$
has an extension to the equivalence relation $\cal H$
in $Q$. Then $M:=Q/\cal H$ is the desired quaternion manifold
with $At (M)=\{ (\Psi _a,U_a): a\in \Lambda  \} $
such that $\Psi _b\circ \Psi _a^{-1}=\Psi _{b,a}$ for each
$U_a\cap U_b\ne \emptyset $, $\Psi _a^{-1}|_{\psi _a(V_a)e}=
\psi _a^{-1}|_{\psi _a(V_a)}$ for each $a$, $\Psi _a^{-1}: Q_a\to U_a$
is the quaternion homeomorphism.
Moreover, each homeomorphism $\psi _a: V_a\to \psi _a(V_a)\subset
\bf C^n$ has the quaternion extension up to the homeomorphism
$\Psi _a: U_a\to \Psi _a(U_a)\subset \bf H^n$. The family of embeddings
$\eta _a: \psi _a(V_a)\hookrightarrow Q_a$ such that
$\pi _{1,1}\circ \eta _a=id$ together with $At(M)$ induces the
complex holomorphic embedding $\theta : N\hookrightarrow M$.
\par {\bf 4.10. Definition.} Let $M$ be a quaternion manifold.
Suppose that for each chart $(U_a,\phi _a)$ of $At(M)$
there exists a quaternion superdifferentiable mapping $\Gamma :
u\in \phi _a(U_a)\mapsto \Gamma (u)\in L_q(X,X,X^*_q;{\bf H})=
L_q(X,X;X)$, where $L_q(X^n,(X^*_q)^m;Y)$ denotes the space
of all quasi-linear mappings from $X^n\times (X^*_q)^m$ into $Y$
(that is, additive and $\bf R$-homogeneous by each argument
$x$ in $X$ or in $X^*_q$), where $X$ and $Y$ are Banach spaces
over $\bf H$, $X^*_q$ denotes the space of all quasi-linear functionals
on $X$ with values in $\bf H$ (see \S 4.1),
$X^*_q=L_q(X;{\bf H})$. If $U_a\cap U_b\ne \emptyset $, let
\par $(1)\quad D(\phi _b\circ \phi _a^{-1}).\Gamma (\phi _a)=
D^2(\phi _b\circ \phi _a^{-1})+\Gamma (\phi _b)\circ (D(\phi _b\circ
\phi _a^{-1})\times D(\phi _b\circ \phi _a^{-1})).$
These $\Gamma (\phi _a)$ are called the Christoffel symbols.
Let ${\cal B}={\cal B}(M)$ be a family of all quaternion holomorphic
vector fields on $M$. For $M$ supplied with $\{ \Gamma (\phi _a): a \} $
define a covariant derivation $(X,Y)\in {\cal B}^2\mapsto \nabla _XY\in
\cal B$:
\par $(2) \nabla _XY(u)=DY(u).X(u)+\Gamma (u)(X(u),Y(u))$,
where $X(u)$ and $Y(u)$ are the principal parts of $X$ and $Y$
on $(U_a,\phi _a)$, $u=\phi _a (z)$, $z\in U_a$.
In this case it is said that $M$ possesses a covariant derivation.
\par {\bf 4.11. Remark.} Certainly for a quaternion manifold
there exists a neighbourhood $V$ of $M$ in $TM$ such that
$\exp : V\to M$ is quaternion holomorphic (see the real case
in \cite{kling}).
\par {\bf 4.12. Theorem.} {\it Let $f$ be a quaternion holomorphic function
such that ${\hat f}$ is quaternion (right) superlinear on a compact
quaternion manifold $M$, then $f$ is constant on $M$.}
\par {\bf Proof.} By the supposition of this theorem
$(f\circ \phi _b^{-1}){\hat .}$ is quaternion
(right) superlinear for each chart $(U_b,\phi _b)$ of $M$.
Since $M$ is compact and $|f(z)|$
is continuous, then there exists a point $q\in M$ at which
$|f(z)|$ attains its maximum. Let $q\in U_b$, then
$f\circ \phi _b^{-1}$ is the quaternion holomorphic function
on $V_b:=\phi _b(U_b)\subset \bf H^n$, where $dim_{\bf H}M=n$.
Consider a polydisk $V$ in $\bf H^n$ with the centre $y=\phi _b(q)$
such that $V\subset V_b$. Put $g(w)=f\circ \phi _b^{-1} (y+(z-y)w)$,
where $w$ is the quaternion variable. Then for each $z\in V$
there exists $\epsilon _z>0$ such that the function $g(w)$
is quaternion holomorphic on the set $W_z:=\{ w: w\in {\bf H},
|w|<1+\epsilon _z \} $
and $|g(w)|$ attains its maximum at $w=0$. In view of Theorem $3.15$
and Remark $3.16$ \cite{luoyst} $g$ is constant on $W_z$, hence
$f$ is constant on $U_b$. By the quaternion holomorphic continuation
$f$ is constant on $M$.

\par Address: Theoretical Department, \\
Institute of General Physics, \\
Russian Academy of Sciences, \\
Str. Vavilov 38, Moscow 119 991 GSP-1, Russia \\
e-mail: ludkovsk@fpl.gpi.ru
\end{document}